\documentclass[11pt]{amsart}

\usepackage[utf8]{inputenc}
\usepackage{graphicx}%
\usepackage{multirow}%
\usepackage{amsmath,amssymb,amsfonts}%
\usepackage{amsthm}
\usepackage{mathrsfs}%
\usepackage[title]{appendix}%
\usepackage{xcolor}%
\usepackage{textcomp}%
\usepackage{manyfoot}%
\usepackage{booktabs}%
\usepackage{listings}
\usepackage{amssymb}
\usepackage{amsfonts}
\usepackage{amsmath}
\usepackage{epsfig}
\usepackage{indentfirst}
\usepackage[usenames,dvipsnames]{pstricks}
\usepackage{pst-grad}
\usepackage{pst-plot}
\usepackage[margin=1.0in]{geometry}
\usepackage{graphicx}
\usepackage{tikz}
\usepackage{tikz-cd}
\usepackage{float}
\usepackage{tabularx}
\usepackage{mathtools}
\theoremstyle{thmstyleone}%
\newtheorem{theorem}{Theorem}[section]
\newtheorem{proposition}[theorem]{Proposition}%
\newtheorem{lemma}[theorem]{Lemma}
\newtheorem{corollary}[theorem]{Corollary}
\newtheorem{remark}[theorem]{Remark}

\theoremstyle{thmstyletwo}%

\theoremstyle{thmstylethree}%
\newtheorem{definition}[theorem]{Definition}%

\DeclareMathOperator {\Hom}{Hom}
\DeclareMathOperator {\End}{End}
\DeclareMathOperator {\Ind}{Ind}
\DeclareMathOperator {\ind}{ind}
\DeclareMathOperator {\Res}{Res}

\DeclareMathOperator {\rad}{rad}
\DeclareMathOperator {\Ext}{Ext}
\DeclareMathOperator {\Stab}{Stab}
\DeclareMathOperator {\Gl}{GL}
\DeclareMathOperator {\arrow}{arrow}

\newcommand{\B}{\mathcal{B}}

\newcommand{\Sg}{\mathfrak{S}}

\usetikzlibrary{%
  matrix,%
  calc,%
  arrows%
}

\title[Permutation modules of the walled Brauer algebras]{Permutation modules of the walled Brauer algebras}
\author{Sulakhana Chowdhury}\address{Indian Institute of Science Education and Research Thiruvananthapuram, Thiruvananthapuram, Email: sulakhana17@iisertvm.ac.in}
\author{Geetha Thangavelu}\address{Indian Institute of Science Education and Research Thiruvananthapuram, Thiruvananthapuram, Email: tgeetha@iisertvm.ac.in }


\begin{document}
\maketitle



\begin{abstract}
In this article, we study the permutation modules and Young modules of the group algebras of the direct product of symmetric groups $K\mathfrak{S}_{a,b}$, and the walled Brauer algebras $\B_{r,t}(\delta)$. In the category of dual Specht-filtered modules, if the characteristic of the field is neither $2$ nor $3$, then the permutation modules are dual Specht filtered, and the Young modules are relative projective cover of the dual Specht modules. We prove that the restriction of the cell modules of $\B_{r,t}(\delta)$ to the group algebras of the direct product of the symmetric groups is dual Specht filtered, and the Young modules act as the relative projective cover of the cell modules of $\B_{r,t}(\delta)$. Finally, we prove that if $\mathrm{char}~K \neq 2,3$, then the permutation module of $\B_{r,t}(\delta)$ can be written as a direct sum of indecomposable Young modules. 
\end{abstract}

\section{Introduction}\label{sec1}
The representation theory of the symmetric groups involves the study of permutation modules, Young modules, and the relations between them. Understanding the decomposition matrix and its properties has yielded profound insights into the representation theory of symmetric groups. It continues to be an active area of research in modern representation theory.

Let $K$ be a field of arbitrary characteristic unless otherwise mentioned. Over fields of positive characteristic, the structure of permutation modules of the group algebras of the symmetric groups $K\Sg_a$ is particularly intriguing. James \cite{JaB}, Klyachko \cite{Kl}, and Grabmeier \cite{Grb} have made significant contributions to the study of the structure of permutation modules of $\Sg_a$ in terms of the Young modules and the Specht modules. Erdmann has also described the Young modules using the representation theory of the symmetric groups \cite{Er}. In this article, we extend this construction to obtain the decomposition of the permutation modules of $K\Sg_{a,b}$ into a direct sum of indecomposables, the Young modules.

Schur-Weyl duality is a fundamental result in the representation theory of Lie groups and diagram algebras, which provides a connection between the representation theory of the symmetric groups and the polynomial representations of the general linear group $\Gl_n$. The Brauer algebras $\mathcal{B}_r(\delta), \delta \in \mathbb{C}$ was first introduced by Richard Brauer in $1937$, restricting the action of $\Gl_n$ in the classical Schur-Weyl duality to orthogonal groups and symplectic groups. Since then, various generalizations of Brauer algebras have been studied relating to the representation theory of groups, such as hyper-octahedral groups, generalized symmetric groups, wreath products of finite groups with symmetric groups, and quantum orthogonal groups. These generalizations include signed Brauer algebras, cyclotomic Brauer algebras, A-Brauer algebras, walled Brauer algebras, and quantum versions of Brauer algebras. One particular case of this duality is the walled Brauer algebras  $\mathcal{B}_{r,t}(\delta), \delta \in \mathbb{C}$, where the general linear group acts on the mixed tensor space $V^{\otimes r} \otimes (V^{*})^{\otimes t}$. The walled Brauer algebras were first introduced independently by Koike \cite{K} and Turaev \cite{T} and have been studied extensively in \cite{BCHL}, \cite{BS}, and by many others. Similar to the Brauer algebras, when the parameter is sufficiently large, the walled Brauer algebras are isomorphic to their corresponding centralizers. They are semisimple except for a few parameters for more details, see \cite{BS}. However, these algebras are known to be non-semisimple for certain small parameters.

Recently, there has been renewed interest in the decomposition of a permutation module of diagram algebras into a direct sum of Young modules. This work was first initiated by Hartmann and Paget in \cite{HP}, in the case of Brauer algebras, and a general construction of permutation modules for all cellularly stratified algebras was given in \cite{HHKP}. Later on, Paul \cite{In} generalized this construction to all cellularly stratified algebras where the input algebras are isomorphic to the group algebras of the symmetric groups or its Hecke algebras using the techniques from \cite{HHKP}. Motivated by these results, we construct the permutation modules of the walled Brauer algebras and obtain its decomposition. Our main result of this article is the following:
\begin{theorem}\label{Dec Per}
If $(l,(\lambda,\mu)) \in \Lambda$ and $\mathrm{char}~K \neq 2,3$, then the decomposition of the permutation module $M(l,(\lambda,\mu))$ is given by $$M(l,(\lambda,\mu)) \cong Y(l,(\lambda,\mu)) \oplus \big ( \bigoplus_{(m,(\lambda',\mu'))}  Y(m,(\lambda',\mu'))^{ \oplus a_{(m,(\lambda',\mu'))}} \big )$$ where $(m,(\lambda',\mu'))$ runs over all $(m,(\lambda',\mu')) \leq  (l,(\lambda,\mu))$, $a_{(m,(\lambda',\mu'))}:=[M(l,(\lambda,\mu)):Y(m,(\lambda',\mu'))]$ is the multiplicity of $Y(m,(\lambda',\mu'))$ in $M(l,(\lambda,\mu))$, and $Y(l,(\lambda,\mu))$ appears exactly once.
\end{theorem}


This article is organized as follows: In Section \ref{Walled Brauer algebras}, we discuss some basic notions and cellularly stratified structure of $\B_{r,t}(\delta)$. In Section \ref{permutation module for direct product group}, we provide a construction for permutation modules and Young modules of $K\mathfrak{S}_{a,b}$. We show that the permutation modules and the Young modules of $K\Sg_{a,b}$ have a dual Specht filtration, and the Young modules are the relative projective covers of the dual Specht modules in $\mathcal{F}_{K\Sg_{r-l,t-l}}(S)$. In Section \ref{Permutation module for Walled Brauer algebra}, we define the permutation modules and the Young modules of $\B_{r,t}(\delta)$, and prove the existence and uniqueness of such Young modules. Moreover, we obtain specific conditions and parameters under which the Young modules can appear in the decomposition of the permutation modules. In Section \ref{imp sec}, we prove Theorem \ref{Dec Per} by obtaining a cell filtration of the permutation modules of $\B_{r,t}(\delta)$ when the characteristic of the field is neither $2$ nor $3$. We also show that if $\mathrm{char}~K \neq 2,3$, then the permutation modules of $\B_{r,t}(\delta)$ and all of its direct summands, in particular, the Young modules admit a cell filtration and are relative projective in $\mathcal{F}_B(\Theta)$. 
%
%
%
%

\section{Walled Brauer algebras} \label{Walled Brauer algebras}
In this section, we recall the definition of the walled Brauer algebras from \cite{CVDM},\cite{RS1}, and prove that these algebras are cellularly stratified.

Let $K$ be a field, and $r$ and $t$ be positive integers, with $\delta$ being a fixed element of $K$. A \textit{$(r,t)$-walled Brauer diagram} is defined as a Brauer diagram consisting of $r+t$ vertices arranged in two rows. The vertices in each row are labeled from left to right as $1,2,\dots, r$ in the top row, and $ r+1,\dots,r+t$ in the bottom row. Edges connecting the vertices within the same row are called ``horizontal edges". We call an edge ``vertical edge" if it connects a vertex in the top row with a vertex in the bottom row. In these diagrams, a wall is placed between the vertices labeled by $r$ and $r+1$, preventing vertical edges from crossing the wall. However, all horizontal edges must cross the wall. 

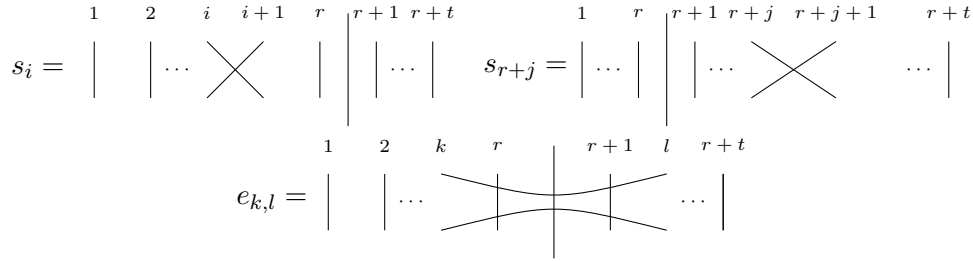
\begin{figure}[H]
\begin{center}
\begin{tikzpicture}[x=0.75cm,y=0.75cm]
\draw[-] (0,0)-- (0,-1);
\draw[-] (1,0)-- (1,-1);
\draw[-] (2,0)-- (3,-1);
\draw[-] (3,0)-- (2,-1);
\draw[-] (4,0)-- (4,-1);
\draw[-] (4.5,0.5)-- (4.5,-1.5);
\draw[-] (5,0)-- (5,-1);
\draw[-] (6,0)-- (6,-1);
\node at (-1,-0.5) {$s_i=$};
\node at (0,0.5) {\tiny $1$};
\node at (1,0.5) { \tiny $2$};
\node at (1.5,-0.5) {\tiny  $\cdots$};
\node at (2,0.5) {\tiny  $i$};
\node at (3,0.5) { \tiny $i+1$};
\node at (4,0.5) {\tiny  $r$};
\node at (5,0.5) { \tiny $r+1$};
\node at (5.5,-0.5) {\tiny  $\cdots$};
\node at (6,0.5) { \tiny $r+t$};
\end{tikzpicture}
\begin{tikzpicture}[x=0.75cm,y=0.75cm]
\draw[-] (0,0)-- (0,-1);
\draw[-] (1,0)-- (1,-1);
\draw[-] (1.5,0.5)-- (1.5,-1.5);
\draw[-] (2,0)-- (2,-1);
\draw[-] (3,0)-- (4.5,-1);
\draw[-] (4.5,0)-- (3,-1);
\draw[-] (6.5,0)-- (6.5,-1);
\node at (-1,-0.5) {$s_{r+j}=$};
\node at (0,0.5) {\tiny $1$};
\node at (0.5,-0.5) {\tiny $\cdots$};
\node at (1,0.5) {\tiny $r$};
\node at (2,0.5) {\tiny $r+1$};
\node at (2.5,-0.5) {\tiny $\cdots$};
\node at (3,0.5) {\tiny $r+j$};
\node at (4.5,0.5) {\tiny $r+j+1$};
\node at (6,-0.5) {\tiny $\cdots$};
\node at (6.5,0.5) {\tiny $r+t$};
\end{tikzpicture}
\begin{tikzpicture}[x=0.75cm,y=0.75cm]
\draw[-] (0,0)-- (0,-1);
\draw[-] (1,0)-- (1,-1);
\draw[-] (2,0).. controls(4,-0.5)..(6,0);
\draw[-] (3,0)-- (3,-1);
\draw[-] (5,0)-- (5,-1);
\draw[-] (4,0.5)-- (4,-1.5);
\draw[-] (7,0)-- (7,-1);
\draw[-] (2,-1).. controls(4,-0.5)..(6,-1);
\draw[-] (7,0)-- (7,-1);
\node at (-1,-0.5) {$e_{k,l}=$};
\node at (0,0.5) {\tiny $1$};
\node at (1,0.5) {\tiny $2$};
\node at (2,0.5) {\tiny $k$};
\node at (1.5,-0.5) {\tiny $\cdots$};
\node at (3,0.5) {\tiny $r$};
\node at (5,0.5) {\tiny $r+1$};
\node at (6.5,-0.5) {\tiny $\cdots$};
\node at (6,0.5) {\tiny $l$};
\node at (7,0.5) {\tiny $r+t$};
\end{tikzpicture}
\caption{Generators of the walled Brauer algebra.}\label{generators of walled Brauer algebra}
\end{center}
\end{figure}

 The \textit{walled Brauer algebra}, denoted by $\mathcal{B}_{r,t} (\delta)$, is an associative $K$-algebra spanned by all $(r,t)$-walled Brauer diagrams. The generators of this algebra are $e_{r,r+1}$, $s_i$ for $1\leq i \leq r+t-1$ with $i\neq r$, and they satisfy the relations given in [\cite{RS1}, Theorem 2.4]. Figure \ref{generators of walled Brauer algebra} provides a visual representation of the generators $s_i$, $s_{r+j}$, and $e_{k,l}$ where $1 \leq i \leq r-1$, $1 \leq j \leq t-1$, and $1\leq k \leq r$, $r+1 \leq l \leq r+t$.

The product of two walled Brauer diagrams $x$ and $y$ is the same as the product of the Brauer diagrams, which is defined as follows: First, place $x$ on top of $y$ and concatenate the bottom vertices of $x$ with the top vertices of $y$. Then remove any closed loops from the concatenated diagram and let $z$ be the resulting diagram. Let $c$ be the number of closed loops removed during this process, we then define $x.y=\delta^{c}z$.

Note that if $d_1$ and $d_2$ are two walled Brauer diagrams having $l_1$ and $l_2$ horizontal edges, respectively, then the product of $d_1$ and $d_2$ will have $l$ horizontal edges, where $l \geq \max \{l_1,l_2\}$. Throughout this article, we call an $K$-linear anti-automorphism $i$ of an algebra $A$ with $i^2=\mathrm{id}$ an \textit{involution}.

\subsection{Cellularly stratification of $\mathcal{B}_{r,t}(\delta)$}

 We recall the definition of cellularly stratified algebras from \cite{HHKP} and prove that walled Brauer algebras are cellularly stratified. 
 
 Let $A, B $ be cellular algebras. We call an algebra $A$ is an inflation of $C$ along $B$, if $A= B\oplus C$ and ideals in the given cell chain must also be ideals in $C$. Also the sections of the chain must satisfy the defining conditions of a cell ideal. For more details reader can refer to [\cite{KX99}, Section 3.3].
\begin{definition} [\cite{HHKP}, Definition 2.1] \label{def cell stra}
Let $A$ be a finite dimensional associative algebra. The algebra $A$ is called \textit{cellularly stratified} if the stratification data $\{B_l,V_l\}_{l=1}^{r}$ satisfies the following:
\begin{itemize}
   \item[(i)] Let $B$ be an inflated algebra (possibly without unit) and $C$ be an algebra (with unit). Define an algebra structure on the vector space $A=B\oplus C$ which extends the two given structures in a way so that $B$ is a two-sided ideal and $A/B\cong C$. We need that $C$ is an ideal, the multiplication is associative, and there exists a unit element of $A$ which maps onto the unit of the quotient $C$. One can see the necessary conditions in \cite{KX99}.
   An inductive application of this procedure on the algebras $C, B_1, B_2, \cdots$ ensures that the layers of the inflations, $B_i= V_i \otimes V_i \otimes B_i^{'}$, are subquotients of ideals in the algebra $A$. We then have that $A$ has an inflation decomposition: $$ A=\bigoplus_{l=1}^{r} V_l \otimes V_l \otimes B_l .$$ 
   \item[(ii)] For each $l=1,2,..., r$, there exist non-zero elements $u_l,v_l \in V_l$ such that $$e_l:= u_l\otimes v_l\otimes \mathit{1}_{B_l}$$ is an idempotent in $A$, where $\mathit{1}_{B_l}$ is an identity element in $B_l$.
   \item[(iii)] For $l>m$, we have $e_l e_m = e_m = e_m e_l$.    
\end{itemize}
\end{definition}
For notational simplicity, we write $s= \min \{r,t\}$, and denote $\mathcal{B}_{r,t} (\delta)=B$. For each $l \in \{1, \cdots, s\}$, we construct an idempotent $e_l$ as follows: when $\delta \neq 0$, the idempotent $e_l$ has $l$ horizontal edges connecting the vertices $r-l+1$ to $r+l$, $r-l+2$ to $r+l-1$, and so on up to $r$ to $r+1$, with the multiplication of factor $\frac{1}{\delta^{l}}$ for $0\leq l \leq s$ (see Figure \ref{idempotents of walled Brauer algebra when delta non zero}). 
\begin{figure}[H]
\begin{center}
\begin{tikzpicture}[x=1cm,y=0.5cm]
\draw[-] (0,0)-- (0,-1);
\draw[-] (1,0)-- (1,-1);
\draw[-] (2,0).. controls(4,-0.5)..(6,0);
\draw[-] (3.5,0).. controls(4,-0.25)..(4.5,0);
\draw[-] (4,0.5)-- (4,-1.5);
\draw[-] (7,0)-- (7,-1);
\draw[-] (2,-1).. controls(4,-0.5)..(6,-1);
\draw[-] (3.5,-1).. controls(4,-0.75)..(4.5,-1);
\draw[-] (7,0)-- (7,-1);
\node at (-1,-0.5) {$e_l=\frac{1}{\delta^{l}}$};
\node at (0,0.5) {\tiny $1$};
\node at (1,0.5) {\tiny $2$};
\node at (2,0.5) {\tiny $r-l+1$};
\node at (2.75,0) {\tiny $\cdots$};
\node at (2.75,-1) {\tiny $\cdots$};
\node at (1.5,-0.5) {\tiny $\cdots$};
\node at (3.5,0.5) {\tiny $r$};
\node at (4.5,0.5) {\tiny $r+1$};
\node at (6.5,-0.5) {\tiny $\cdots$};
\node at (5.5,0) {\tiny $\cdots$};
\node at (5.5,-1) {\tiny $\cdots$};
\node at (6,0.5) {\tiny $r+l$};
\node at (7,0.5) {\tiny $r+t$};
\end{tikzpicture}
\caption{Idempotents of $\mathcal{B}_{r,t}(\delta)$ for $\delta \neq 0$.}
\label{idempotents of walled Brauer algebra when delta non zero}
\end{center}
\end{figure}
However, when $\delta =0$ and one of $r$ or $t$ is at least $2$, the idempotent $e_l$ is defined as shown in Figure \ref{idempotents of walled Brauer algebra when delta zero}.
\begin{figure}[H]
\begin{center}
\begin{tikzpicture}[x=0.9cm,y=0.5cm]
\draw[-] (0,0)-- (0,-1);
\draw[-] (1,0)-- (1,-1);
\draw[-] (2,0).. controls(4,-0.5)..(6,0);
\draw[-] (2,-1).. controls(4,-0.65)..(7,-1);
\draw[-] (3.5,0).. controls(4,-0.25)..(4.5,0);
\draw[-] (4,0.5)-- (4,-1.5);
\draw[-] (7,0)-- (6,-1);
\draw[-] (8,0)-- (8,-1);
\draw[-] (3.5,-1).. controls(4,-0.85)..(4.5,-1);
\node at (-1,-0.5) {$e_l=$};
\node at (0,0.5) {\tiny $1$};
\node at (1,0.5) {\tiny $2$};
\node at (1.5,-0.5) {\tiny $\cdots$};
\node at (2,0.5) {\tiny $r-l+1$};
\node at (2.75,0) {\tiny $\cdots$};
\node at (3.5,0.5) {\tiny $r$};
\node at (4.5,0.5) {\tiny $r+1$};
\node at (5.25,0) {\tiny $\cdots$};
\node at (6,0.5) {\tiny $r+l$};
\node at (7.20,-0.5) {\tiny $\cdots$};
\node at (7,0.5) {\tiny $r+l+1$};
\node at (8,0.5) {\tiny $r+t$};
\end{tikzpicture}
\caption{Idempotents of $\B_{r,t}(\delta)$ for $\delta= 0$ and when one of $r$ or $t$ is at least $2$.}\label{idempotents of walled Brauer algebra when delta zero}
\end{center}
\end{figure}

Let $V_l$ be the $K$-vector space spanned by all $(r,t)$-partial diagrams that have $l$ horizontal edges with $r-l$ free vertices on the left side of the wall and $t-l$ free vertices on the right side of the wall. For $v \in V_l$ and $d\in \B_{r,t}(\delta)$, the action of $d$ on $v$ is defined by placing $v$ above $d$ and concatenating them  in same manner as in $\B_{r,t}(\delta)$. The product is considered zero if the resulting configuration does not have exactly $l$ horizontal edges. We define an involution $i$ by flipping the diagram along a horizontal line. Now, consider a filtration of $B$ by cell ideals: $$ 0\subseteq J_s \subseteq J_{s-1} \subseteq \cdots \subseteq J_{0} = B,$$ where $J_l$ is a two-sided ideal of $B$ containing at least $l$ horizontal edges. The subquotient $J_l /J_{l+1}$ is isomorphic to the inflation $V_{l}\otimes V_{l} \otimes K\mathfrak{S}_{r-l,t-l}$ along the vector space $V_l$. The following proposition satisfies the first condition of Definition \ref{def cell stra}.

\begin{proposition} [\cite{CVDM}, Proposition 2.6]
The walled Brauer algebra can be written as an iterated inflation of $K\mathfrak{S}_{r-l,t-l}$ along $V_l$ $$ \mathcal{B}_{r,t}(\delta) =\bigoplus_{l=0}^{s} V_{l}\otimes V_{l} \otimes K\mathfrak{S}_{r-l,t-l}, \text{ where } s= \min \{r,t\}.$$ 
\end{proposition} 

The group algebra $K\mathfrak{S}_{a,b}$ is cellular due to the fact that the direct product of cellular algebras is cellular as shown in [\cite{GM}, Example 3.1.6]. The cell ideals are given by $S^{\lambda} \boxtimes_{K} S^{\mu}$, where $\lambda \vdash r-l$ and $\mu \vdash t-l$. The second condition for Definition \ref{def cell stra} is satisfied by constructing the idempotent $e_l$ when $\delta\neq 0$. This condition also holds when $\delta=0$ and one of $r$ or $t$ is at least greater than or equal to $2$. However, when $r=t=1$, the two-sided ideal $J_1$ becomes nilpotent and therefore cannot be generated by the idempotent $e_1$. We establish the third condition of Definition \ref{def cell stra} by utilizing the fact that the number of horizontal edges in the product of two diagrams is greater than or equal to the number of horizontal edges in each diagram. Since the idempotents $e_l$ satisfy this property, the third condition holds. Finally, we conclude that

\begin{theorem}\label{cell stra}
Let $K$ be a field and $\delta \in K$. If $\delta \neq 0$, then $\B_{r,t}(\delta)$ is cellularly stratified. If $\delta=0$, and one of $r$ and $t$ is at least $2$, then $\B_{r,t}(\delta)$ is cellularly stratified. 
\end{theorem} 

\section{Permutation modules of $K\mathfrak{S}_{a,b}$ } \label{permutation module for direct product group}

In this section, we study the decomposition of the permutation modules of $K\mathfrak{S}_{a,b}$ into indecomposable Young modules. We prove this by extending the construction by James in \cite{Ja}. We begin by reviewing the structure of the permutation modules and the Young modules of the symmetric groups. Towards the end of the section, we prove that the permutation modules and the Young modules of $K\Sg_{a,b}$ have a dual Specht filtration and are relative projective in the category of modules with a dual Specht filtration. 

\subsection{Representation of symmetric group $\mathfrak{S}_a$ } 
Let $a$ be a natural number. A \textit{composition} $\lambda$ of $a$ (denoted by $\lambda \vDash a$) is a finite sequence of non-negative integers $\lambda = (\lambda_1, \cdots, \lambda_k)$ such that $\sum_{i=1}^k |\lambda_i |= a$. A \textit{partition} $\lambda$ of a natural number $a$ (denoted by $\lambda \vdash a$) is a composition of $a$ whose parts are strictly positive and arranged in non-increasing order. A \textit{Young diagram} $[\lambda]$ of shape $\lambda$ is an array of boxes $\{(i,j) : i,j \in \mathbb{N}, 1 \leq i \leq k, 1 \leq j \leq \lambda_i\}$. The \textit{conjugate} of $\lambda$, denoted by $\lambda'$, is obtained by interchanging the rows and columns of $\lambda$. The set of all partition of $a$ is denoted by $\Lambda_a$. For $\lambda, \mu \in \Lambda_a$, the dominance order on the index set $\Lambda_a$ is given by $\lambda \unrhd \mu$ if $\sum_{i=1}^l |\lambda_i| \geq \sum_{i=1}^l |\mu_i|$ for all $l$. A $\lambda$-\textit{tableau} $\mathfrak{t}$ is a Young diagram of shape $\lambda$ where the entries are filled with numbers $1,2, \cdots, a$ without repetition. The symmetric group $\mathfrak{S}_a$ acts from the right on a $\lambda$-tableau by permuting the entries. Given a $\lambda$-tableau $\mathfrak{t}$, the \textit{row stabilizer}  $R_{\mathfrak{t}}$ of  $\mathfrak{t}$ is a subgroup of $\mathfrak{S}_a$ that fixes rows of $\mathfrak{t}$ as a set, and the \textit{column stabilizer} $C_{\mathfrak{t}}$ can be defined similarly. Let $\mathfrak{t}_1, \mathfrak{t}_2$ be two $\lambda$-tableaux. An equivalence relation ``$\sim$" on $\lambda$-tableaux is defined by: $\mathfrak{t}_1 \sim \mathfrak{t}_2$ if and only if there exists $\sigma \in R_{\mathfrak{t}}$ such that $\mathfrak{t}_2 \cdot \sigma = \mathfrak{t}_1$. We call the equivalence class that contains $\mathfrak{t}$ to be the \textit{$\lambda$-tabloid}, denoted by $\{\mathfrak{t}\}$. The symmetric group $\Sg_a$ acts on the set of $\lambda$-tabloids by $\{\mathfrak{t}\}\cdot \sigma=\{\mathfrak{t}\cdot \sigma\}$. For each $\lambda \vdash a$, the \textit{Young permutation module} $M^{\lambda}$ is a right $K\Sg_a$-module where the underlying $K$-vector space is spanned by all $\lambda$-tabloids. The \textit{Young subgroup} $\mathfrak{S}_{\lambda}$ of $\mathfrak{S}_a$ is defined as $\mathfrak{S}_{\lambda}:= \mathfrak{S}_{\{1,2,\cdots, \lambda_1\}} \times \mathfrak{S}_{\{\lambda_1+1, \cdots, \lambda_1+\lambda_2\}} \times \cdots \times \mathfrak{S}_{\{\lambda_1+\cdots+ \lambda_{k-1}+1, \cdots, \lambda_k\}}$. Alternatively, one can view this permutation module as $M^{\lambda} \cong \Ind_{K\Sg_{\lambda}}^{K\mathfrak{S}_a} 1_{K\Sg_{\lambda}} =  1_{K\Sg_{\lambda}}\otimes_{K\Sg_{\lambda}}K\mathfrak{S}_a $, where $1_{K\Sg_{\lambda}}$ is the trivial $K\Sg_{a}$-module.

The \textit{alternating column sum} $k_\mathfrak{t}$ is an element of $\mathfrak{S}_a$ obtained by taking the alternating sum of the elements in $C_{\mathfrak{t}}$, i.e., $ k_{\mathfrak{t}}= \sum_{\sigma\in C_\mathfrak{t}} (\mathrm{sgn}~ \sigma)\sigma$. The $\lambda$-\textit{polytabloid} $e_{\mathfrak{t}}$ is an element of $M^{\lambda}$, defined by $e_{\mathfrak{t}}:= \{\mathfrak{t}\}\cdot k_{\mathfrak{t}} $. The \textit{Specht module} $S^{\lambda}$ is a cyclic submodule of $M^{\lambda}$ generated by the $\lambda$-polytabloid $e_{\mathfrak{t}}$. Over a field of characteristic $0$, these modules are self-dual and absolutely irreducible. 
%
 A partition $\lambda=(\lambda_1,\lambda_2, \cdots )$ of $a$ is \textit{$p$-regular} if it has at most $p-1$ parts of any given size.  
If $\mathrm{char}~K=p$ is a prime, then the set $\{D^\lambda=S^{\lambda}/\rad~S^{\lambda}: \lambda \text{ is a }p\text{-regular partition of }a\}$ forms a complete set of non isomorphic simple $K\Sg_a$-modules by [\cite{JaB}, Theorem 11.5 ].

Using [\cite{CR}, Theorem 6.12], the permutation module can be expressed as a direct sum of indecomposable modules, i.e., $M^{\lambda} = \bigoplus Y_i$. By James submodule theorem, there exists a unique indecomposable module $Y_i$ containig the Specht module $S^{\lambda}$. We refer to this module as the \textit{Young module} and denote it by $Y^{\lambda}$. The classification of the Young modules of $K\Sg_a$ given by James is as follows:

\begin{theorem}[\cite{Ja}, Theorem 3.1]\label{De sym}
Let $\mathrm{char}~K=p$ and $\lambda, \mu \vdash n$. The permutation module $M^{\lambda}$ of $\Sg_a$ can be decomposed as a direct sum of Young modules $\{Y^{\mu}: \mu \unrhd \lambda \}$, that is $$M^{\lambda}\cong Y^{\lambda} \oplus \big(\bigoplus_{\mu \unrhd \lambda} Y^{\mu^{ \oplus a_\mu}}\big),$$ where 
\begin{itemize}
\item[1.] $Y^{\lambda}$ appears exactly once and the multiplicity of $Y^{\mu}$ in $M^{\lambda}$, $[M^{\lambda}: Y^{\mu}]$ is called the $p$-Kostka number, and is denoted by $a_{\mu}$. 
\item[2.] Young modules are pairwise inequivalent i.e. $Y^{\lambda}\cong Y^{\mu}$ if and only if $\lambda=\mu$.
\end{itemize}
\end{theorem}
Note that, over the field of characteristic $0$, the Young modules and the Specht modules coincide. 

Furthermore, the \textit{dual Specht module}, $S_{\lambda}$, can be written as $ S_{\lambda}:= (S^{\lambda})^*\cong  S^{\lambda'}\otimes \mathrm{ sgn}_a$ as shown in [\cite{JaB}, Theorem 8.15], where $\lambda'$ is the conjugate of the partition $\lambda$. We say that a $K\Sg_a$-module $M$ has a \textit{dual Specht filtration} if $M$ has a chain of submodules: 
\[
M=M_0 \supseteq M_1 \supseteq \cdots \supseteq M_m \supseteq M_{m+1}=\{0\}
\]
with each subquotient $M_i/M_{i+1}$ isomorphic to some dual Specht module of $K\Sg_a$. Consider $\mathcal{F}_{K\Sg_a}(S)$ to be the category of $K\Sg_a$-modules with a dual Specht filtration.

\begin{lemma}\label{sym grp has dual Specht filtration} If $\mathrm{char}~K \neq 2,3$, then the permutation module of $K\mathfrak{S}_a$ has a dual Specht filtration. In particular, $M^{\lambda} \in \mathcal{F}_{K\Sg_a}(S)$.
\end{lemma}
\begin{proof}
Given any $\alpha \vdash a$ and $\lambda \vDash a$, by [\cite{HN}, Proposition 4.2.1], we have $$\Ext_{K\mathfrak{S}_a}^{1}(S^{\alpha}, M^{\lambda})=0.$$ Thus, the result follows from [\cite{H}, Theorem 3.4,(ii)].
\end{proof}

\subsection{Representation theory of $K\mathfrak{S}_{a,b}$} \label{prop dipdt}

In this section, we will describe the analog of Theorem \ref{De sym} for $K\mathfrak{S}_{a,b}$. Let $a,b \in \mathbb{N}$, and define $K\Sg_{a,b}$ to be the group algebra of the direct product of two symmetric groups $\mathfrak{S}_a$ and $\mathfrak{S}_b$, consisting of the elements of the form  $(\sigma, \tau )$, where $ \sigma \in \mathfrak{S}_a$ and $ \tau \in \mathfrak{S}_b $. 

For $\lambda \vdash a$ and $\mu \vdash b$, we define the index set $\Lambda_{a,b}$ as $ \Lambda_{a,b}:= \{ (\lambda,\mu) : \lambda \vdash a, \mu \vdash b \}$. This is a subset of the set of all \textit{bi-partitions} of $a+b$, denoted by $\Lambda_{a+b}$. Let $(\lambda,\mu), (\lambda',\mu') \in \Lambda_{a,b}$, the \textit{dominance order} $``\unrhd"$ on $ \Lambda_{a,b}$ is defined by $$ (\lambda,\mu) \unrhd (\lambda',\mu') \text{ if } \lambda \unrhd \lambda' \text{ and } \mu \unrhd \mu'.$$ Given the tabloids $\{t_{\lambda}\} \in M^{\lambda}$ and $\{t_{\mu}\} \in M^{\mu}$, define a \textit{$(\lambda,\mu)$-tabloid} as a pair $(\{t_{\lambda}\},\{t_{\mu}\})$, or equivalently a horizontal juxtaposition of $\{t_{\lambda}\}$ and $\{t_{\mu}\}$. Let $(\lambda,\mu) \in \Lambda_{a,b}$, and consider $M_{\lambda, \mu}$ to be a $K$-vector space spanned by outer tensor product of all $\lambda$ and $\mu$-tabloids, i.e., $M^{\lambda,\mu}=\{ \{\mathfrak{t}_{\lambda,\mu}\}:= (\{\mathfrak{t}_{\lambda}\}\boxtimes_{K} \{\mathfrak{t}_{\mu}\}): (\lambda,\mu) \in \Lambda_{a,b} \}.$ The group algebra $K\Sg_{a,b}$ acts naturally on $M^{\lambda,\mu}$ by $$\{\mathfrak{t}_{\lambda,\mu}\} \cdot (\sigma,\tau)= (\{\mathfrak{t}_{\lambda}\}\boxtimes_{K} \{\mathfrak{t}_{\mu}\})\cdot (\sigma,\tau)= (\{\mathfrak{t}_{\lambda}\cdot \sigma\}\boxtimes_{K} \{\mathfrak{t}_{\mu}\cdot \tau\}).$$ Thus, $M^{\lambda,\mu}$ is a $K\Sg_{a,b}$-module, which is the outer tensor product of two permutation modules $M^{\lambda}$ and $M^{\mu}$ of $K\Sg_{a}$ and $K\Sg_b$, respectively. The $K\Sg_{a,b}$-module $M^{\lambda,\mu}$ is called a \textit{Young permutation module}. 

Given $(\lambda,\mu) \in \Lambda_{a,b}$, we denote $\Sg_{\lambda,\mu}$ to be the subgroup of $\Sg_{a,b}$ defined by $\Sg_{\lambda,\mu}= \Sg_{\lambda} \times \Sg_\mu$. The subgroups $\Sg_{\lambda,\mu}$ are called the \textit{Young subgroup}s of $\Sg_{a,b}$ for all $(\lambda,\mu) \in \Lambda_{a,b}$. 

\begin{lemma}
 If $(\lambda,\mu) \in \Lambda_{a,b}$, then we have an isomorphism of $K\Sg_{a,b}$-modules $$ M^{\lambda,\mu} \cong \Ind_{K\Sg_{\lambda,\mu}}^{K\mathfrak{S}_{a,b}} 1_{K\Sg_{\lambda,\mu}}.$$
\end{lemma}
\begin{proof}
The proof can be obtained as follows:
\begin{align*}
M^{\lambda,\mu}=M^{\lambda} \boxtimes_K M^{\mu} &\cong (\Ind_{K\Sg_{\lambda}}^{K\mathfrak{S}_a} 1_{K\Sg_\lambda}) \boxtimes_{K} (\Ind_{K\Sg_{\mu}}^{K\mathfrak{S}_b} 1_{K\Sg_\mu})\\
& \cong \Ind_{K\Sg_{\lambda}\times K\Sg_{\mu}}^{K\mathfrak{S}_a \times K\mathfrak{S}_b} 1_{K\Sg_{\lambda}\times K\Sg_{\mu}}\\
&=\Ind_{K\Sg_{\lambda,\mu}}^{K\mathfrak{S}_{a,b}} 1_{K\Sg_{\lambda,\mu}},
\end{align*}
where the second congruence holds because the outer tensor product commutes with induction as shown in [\cite{CR}, Lemma 10.17].
\end{proof}

Let $k_{\mathfrak{t}_{\lambda}}$ and $k_{\mathfrak{t}_{\mu}}$ represent the alternating column sum of the $\lambda$-tableau $\mathfrak{t}_{\lambda}$ and the $\mu$-tableau $\mathfrak{t}_{\mu}$, respectively. The \textit{alternating column sum} of the $(\lambda,\mu)$-tableau $\mathfrak{t}_{\lambda,\mu}$ is an element of $K\Sg_{a,b}$, and is defined by $k_{\mathfrak{t}_{\lambda,\mu}}= \sum_{\sigma \in C_{\mathfrak{t}_{\lambda}}, \tau \in C_{\mathfrak{t}_{\mu}}} (\text{sgn } \sigma\cdot \text{sgn } \tau) (\sigma,\tau)$, where $C_{\mathfrak{t}_{\lambda}}$ and $C_{\mathfrak{t}_{\mu}}$ are the column stabilizers of the tableaux $\mathfrak{t}_{\lambda}$ and $\mathfrak{t}_{\mu}$, respectively. The \textit{$(\lambda,\mu)$-polytabloid} $e_{\mathfrak{t}_{\lambda,\mu}}$ is defined by the product $ \{\mathfrak{t}_{\lambda,\mu}\}\cdot k_{\mathfrak{t}_{\lambda,\mu}}$. Now, consider the $K\Sg_{a,b}$-submodule $S^{\lambda,\mu}$ of $M^{\lambda,\mu}$ generated by all the polytabloids $e_{\mathfrak{t}_{\lambda,\mu}}$, and we call $S^{\lambda,\mu}$ the \textit{Specht module} of $K\mathfrak{S}_{a,b}$. 

\begin{lemma}
The Specht module $S^{\lambda, \mu}$ of $K\Sg_{a,b}$ is the outer tensor product of $S^{\lambda}$ and $S^{\mu}$, and it is a cyclic module.
\end{lemma}
\begin{proof}
Let $e_{\mathfrak{t}_{\lambda}}, e_{\mathfrak{t}_{\mu}}$ and $e_{\mathfrak{t}_{\lambda,\mu}}$ be the generator of $S^{\lambda}$,  $S^{\mu}$ and $S^{\lambda,\mu}$,respectively. Then we have 
\begin{align*}
e_{\mathfrak{t}_{\lambda,\mu}}&= \{\mathfrak{t}_{\lambda,\mu}\}\cdot k_{\mathfrak{t}_{\lambda,\mu}}\\
&=\sum_{\sigma \in C_{\mathfrak{t}_{\lambda}}, \tau \in C_{\mathfrak{t}_{\mu}}} (\{\mathfrak{t}_{\lambda}\} \boxtimes_K \{ \mathfrak{t}_{\mu}\})\cdot (\text{sgn } \sigma\cdot \text{sgn } \tau) (\sigma,\tau)\\
&= \big( \sum_{\sigma \in C_{\mathfrak{t}_{\lambda}}} \{\mathfrak{t}_{\lambda}\}\cdot (\text{sgn } \sigma) \sigma \big) \boxtimes_K  \big (\sum_{\tau \in C_{\mathfrak{t}_{\mu}}} \{\mathfrak{t}_{\mu}\}\cdot  (\text{sgn } \tau) \tau \big)\\
&= e_{\mathfrak{t}_{\lambda}} \boxtimes_K e_{\mathfrak{t}_{\mu}}. 
\end{align*}
Therefore, we have $S^{\lambda,\mu} \cong S^{\lambda} \boxtimes_K S^{\mu}$. Since $S^{\lambda}$ and $S^{\mu}$ are cyclic modules, the outer tensor product of two cyclic modules is again cyclic. Thus, $S^{\lambda,\mu}$ is also cyclic.
\end{proof}

\begin{lemma} 
If $K$ has characteristic $0$, then the Specht modules $S^{\lambda,\mu}$, for $(\lambda,\mu) \in \Lambda_{a,b}$ form a complete set of simple $K\mathfrak{S}_{a,b}$-modules.
\end{lemma}
\begin{proof}
When $\mathrm{char}~K=0$, by [\cite{JaB}, Theorem 4.12], the Specht modules $S^{\lambda}$ and $S^{\mu}$ are simple. Then by using [\cite{CR}, Theorem 10.33], we conclude that  $S^{\lambda,\mu}=S^{\lambda} \boxtimes_K S^{\mu}$ is a simple $K\mathfrak{S}_{a,b}$-module, where the Specht modules $S^{\lambda}$ and $S^{\mu}$ are uniquely determined up to isomorphism.
\end{proof}

We can define a bilinear form $\langle-,-\rangle$ on the Young permutation module
$M^{\lambda,\mu}$ by requiring that it is orthonormal on tabloids: That is, $\langle \{\mathfrak{t}_{\lambda,\mu}\}, \{\mathfrak{t}_{\lambda',\mu'}\} \rangle=1 $ if $\{\mathfrak{t}_{\lambda,\mu}\}= \{\mathfrak{t}_{\lambda',\mu'}\}$, i.e., $\{\mathfrak{t}_{\lambda}\}=\{\mathfrak{t}_{\lambda'}\}$ and $\{\mathfrak{t}_{\mu}\}=\{\mathfrak{t}_{\mu'}\}$, otherwise $0$. 

For $(\lambda, \mu), (\lambda',\mu') \in \Lambda_{a,b}$ and for any submodule $N$ of $M^{\lambda,\mu}$, we define the radical of $N$ as: $$\rad~N= \big\{ \{t_{\lambda,\mu}\} \in M^{\lambda,\mu}: \langle \{t_{\lambda,\mu}\}, \{t_{\lambda',\mu'}\} \rangle =0, \text{ for all }\{t_{\lambda',\mu'}\} \in N \big \}.$$
\begin{lemma} 
If $\mathrm{char}~K=p$, a prime, then the modules $D^{\lambda,\mu}:= S^{\lambda,\mu}/\rad~S^{\lambda,\mu}$ form a complete set of simple $K\mathfrak{S}_{a,b}$-modules, where $ \lambda$ and $\mu$ are the $p$-regular partitions of $a$ and $b$, respectively.
\end{lemma} 
\begin{proof}
By applying [\cite{CR}, Theorem 7.10], we know that $K\big(\mathfrak{S}_{a}/ \text{rad} \mathfrak{S}_{a}\big)$ and $K\big(\mathfrak{S}_{b}/ \rad \mathfrak{S}_{b}\big)$ are separable $K$-algebras, since $K\big(\mathfrak{S}_{a}/ \rad \mathfrak{S}_{a}\big)$, and $K\big(\mathfrak{S}_{b}/ \rad \mathfrak{S}_{b}\big)$ are always semi-simple, meaning every field is a splitting field. If $K$ is of characteristic $p$ and for $\lambda$ and $\mu$, $p$-regular partitions of $a$ and $b$ respectively, then $S^{\lambda}/\rad  S^{\lambda}$ and $S^{\mu}/\rad S^{\mu}$ respectively are the simple $K\mathfrak{S}_{a}$ and $K\mathfrak{S}_{b}$-modules respectively. By applying [\cite{CR}, Theorem 10.38], the simple $K\mathfrak{S}_{a,b}$-module is given by: 
 \begin{align*}
(S^{\lambda}/\text{rad } S^{\lambda}) \boxtimes (S^{\mu}/(\text{rad}S^{\mu})
&=(S^{\lambda}\boxtimes S^{\mu})/\big( (\text{rad } S^{\lambda}) \boxtimes (\text{rad}S^{\mu})\big)~[\text{by [\cite{CR}, } 10.39]] \\
&=S^{\lambda,\mu}/ \text{rad}(S^{\lambda} \boxtimes S^{\mu})\\
&=S^{\lambda,\mu}/\text{rad}~ S^{\lambda,\mu}.\nolinebreak
\end{align*}
\end{proof}  

Now, we proceed to the decomposition of the permutation module into a direct sum of indecomposable $K\mathfrak{S}_{a,b}$-modules. Write the permutation module $M^{\lambda,\mu}$ as a direct sum of indecomposable $K\Sg_{a,b}$-modules, say $M^{\lambda,\mu}= \bigoplus_i Y_i$. Let $k_{\mathfrak{t}_{\lambda,\mu}}$ be the alternating column sum. We have
\begin{align*}
M^{\lambda,\mu} \cdot k_{\mathfrak{t}_{\lambda,\mu}}&= (M^{\lambda} \boxtimes M^{ \mu}) \cdot (\sum_{\sigma\in C_{\mathfrak{t}_{\lambda}}, \tau \in C_{\mathfrak{t}_{\mu}}} (\mathrm{sgn}~\sigma \cdot \mathrm{sgn}~\tau)( \sigma,\tau))\\
&= \big( M^{\lambda} \cdot \sum_{\sigma\in C_{\mathfrak{t}_{\lambda}}}(\mathrm{sgn}~ \sigma) \sigma \big) \boxtimes_K \big( M^{\lambda}\cdot \sum_{\tau\in C_{\mathfrak{t}_{\mu}}}(\mathrm{sgn}~ \tau) \tau\big)\\
&=C e_{\mathfrak{t}_{\lambda}} \boxtimes C'e_{\mathfrak{t}_{\mu}}\\
&= K e_{\mathfrak{t}_{\lambda, \mu}}, 
\end{align*}
where $C$ and $C'$ are non-zero constants. Therefore, there exists a unique summand $Y_j$ such that $S^{\lambda,\mu} \cap Y_j \neq 0$ (since $ 0 \neq S^{\lambda,\mu} \cap Y_j = ( S^{\lambda} \cap Y_j^{1} ) \boxtimes (S^{\mu} \cap Y_j^{2})$), where $Y_j^1$ (resp. $Y_j^2$) is the unique summand of $M^{\lambda}$ (resp. $M^{\mu}$) containing $S^{\lambda}$ (resp. $S^{\mu}$) as a submodule, by James submodule theorem [\cite{JaB}, Theorem 4.8].

\begin{theorem}\label{de per dirpdt}
Let $K$ be a field of characteristic $p$. Let $(\lambda,\mu), (\lambda',\mu') \in \Lambda_{a,b}$, the permutation module $M^{\lambda,\mu}$ can be decomposed into a direct sum of indecomposable Young modules $Y^{\lambda,\mu}$, 
\begin{itemize}
\item[1.] The indecomposable summands are of the form $\{Y^{\lambda',\mu'}:= Y^{\lambda'} \boxtimes Y^{\mu'}: \lambda' \unrhd \lambda, \mu' \unrhd \mu\}$.
\item[2.] $Y^{\lambda,\mu}$ will appear exactly once, and the multiplicity of $Y^{\lambda',\mu'}$ in $M^{\lambda,\mu}$ is given by $a_{\lambda',\mu'}:=[M^{\lambda,\mu} : Y^{\lambda',\mu'}]$, which is a non-zero integer.
\item[3.] Any two Young modules $Y^{\lambda',\mu'}$ and $Y^{\lambda,\mu}$ are isomorphic if and only if $\lambda'=\lambda$ and $ \mu'=\mu $.
\end{itemize}
Together we write 
\begin{equation*}
M^{\lambda,\mu} \cong Y^{\lambda,\mu} \oplus \big(\bigoplus_{\substack{\lambda' \unrhd \lambda\\ \mu' \unrhd \mu}}  Y^{\lambda',\mu'^{\oplus a_{\lambda',\mu'}} } \big).
\end{equation*}
\end{theorem}

\begin{proof}

\begin{itemize}
\item[1.] By applying the Theorem \ref{De sym}, we have
 \begin{align*}
M^{\lambda,\mu}= M^{\lambda} \boxtimes M^{\mu} &= (\bigoplus_{\lambda' \unrhd \lambda}Y^{\lambda'^{\oplus a_{\lambda'}}})\boxtimes  (\bigoplus_{\mu' \unrhd \mu} Y^{\mu'^{\oplus a_{\mu'}}} )\\
&=\bigoplus_{\substack{(\lambda',\mu') \unrhd (\lambda, \mu)}} (Y^{\lambda'^{\oplus a_{\lambda'}}} \boxtimes Y^{\mu'^{\oplus a_{\mu'}}})\\
&= \bigoplus_{\substack{(\lambda',\mu') \unrhd (\lambda, \mu)}} (Y^{\lambda'} \boxtimes Y^{\mu'})^{ \oplus a_{\lambda',\mu'}},
\end{align*}
where $a_{\lambda',\mu'}$ is the product of $a_{\lambda'}$ and $a_{\mu'}$, and in the last equality $Y^{\lambda'} \boxtimes Y^{\mu'}$ is indecomposable, which follows from the fact that the outer tensor product of two indecomposable modules is indecomposable. So we conclude that the indecomposable summands are parametrized by $(\lambda,\mu) \in \Lambda_{a,b}$, where $(\lambda',\mu') \unrhd (\lambda, \mu) $.
\item[2.] The multiplicity of each indecomposable summand in the decomposition of the permutation module is given as follows:
\begin{align*}
[M^{\lambda,\mu}: Y^{\lambda,\mu}]&=\text{dim}_{K} (\Hom_{K\mathfrak{S}_{a,b}} (M^{\lambda,\mu}, Y^{\lambda,\mu}))\\
&=\text{dim}_{K} (\Hom_{K\mathfrak{S}_a} (M^{\lambda}, Y^{\lambda}) \boxtimes \Hom_{K\mathfrak{S}_b} (M^{\mu}, Y^{\mu}))\\
&=\text{dim}_{K} (\Hom_{K\mathfrak{S}_a} (M^{\lambda}, Y^{\lambda}))  \text{dim}_{K} (\Hom_{K\mathfrak{S}_b} (M^{\mu}, Y^{\mu}))\\
&=[M^{\lambda}:Y^{\lambda}]\cdot[M^{\mu}:Y^{\mu}]\\
&= 1,
\end{align*}
where second equality holds by [\cite{Xi}, Lemma 3.2], and the last equality is due to part 1 of Theorem \ref{De sym}. Hence, for $\lambda' \unrhd \lambda$ and $ \mu' \unrhd \mu $, the multiplicity of $Y^{\lambda',\mu'}$ in $M^{\lambda,\mu}$ is $a_{\lambda',\mu'}$ which is the product of $a_{\lambda'}$ and $a_{\mu'}$, and is a non-negative integer by using Theorem \ref{De sym} part 1.
\item[3.]
To determine when two Young modules $Y^{\lambda',\mu'}$ and $Y^{\lambda,\mu}$ are isomorphic, assume that $Y^{\lambda',\mu'} \cong Y^{\lambda,\mu}$. Since $Y^{\lambda,\mu}$ is an indecomposable summand of $M^{\lambda,\mu}$, therefore $Y^{\lambda',\mu'}$ must also appear in the decomposition of $M^{\lambda,\mu}$. Since we know that all the indecomposable summands are parametrized by $(\lambda',\mu')$ where $\lambda' \unrhd \lambda, \mu' \unrhd \mu$, so we must have $(\lambda',\mu') \unrhd (\lambda, \mu)$. By symmetry, $(\lambda,\mu) \unrhd (\lambda', \mu')$. Therefore, $\lambda'=\lambda$ and $\mu'=\mu$. Thus, the Young modules $Y^{\lambda',\mu'}$ and $Y^{\lambda,\mu}$ are isomorphic if and only if $\lambda'=\lambda$ and $\mu'=\mu$.
\end{itemize}
\end{proof}
 We call the unique summand $Y^{\lambda}\boxtimes Y^{\mu}$ to be the \textit{Young module} of $K\mathfrak{S}_{a,b}$ and denote it by $Y^{\lambda,\mu}$ analogous to James \cite{Ja}. 

Next, we define the \textit{dual Specht module} $S_{\lambda,\mu}$ in a similar manner. 
\begin{lemma}\label{dual Specht module of direct product is tensor product of dual Specht module for direct product}
For $(\lambda,\mu) \in \Lambda_{a,b}$, the dual Specht module $S_{\lambda,\mu}$ is isomorphic to the outer tensor product of the dual Specht modules $S_{\lambda} $ of $K\mathfrak{S}_a$  and $S_{\mu}$ of $K\mathfrak{S}_b$. 
\end{lemma}
\begin{proof}
It is known that $ S_{\lambda,\mu} := (S^{\lambda,\mu})^{*} \cong (S^{\lambda})^{*} \boxtimes (S^{\mu})^{*} = S_{\lambda} \boxtimes S_{\mu}$. This holds because the outer tensor product of dual Specht modules is isomorphic to a dual of the outer tensor product of Specht modules.
\end{proof}

We end this section by discussing some homological properties of permutation modules and Young modules of $K\Sg_{a,b}$, which will be useful in Section \ref{imp sec}.
\begin{lemma}\label{tensor product of two filtration for walled Brauer}
If $M \in \mathcal{F}_{K\Sg_a}(S)$ and $N \in \mathcal{F}_{K\Sg_b}(S)$, then $M \boxtimes N$ is a $K \Sg_{a,b}$-module, and $M\boxtimes N \in \mathcal{F}_{K\Sg_a}(S) \times \mathcal{F}_{K\Sg_b}(S)$.
\end{lemma}
\begin{proof}
This follows from [\cite{CT}, Lemma 4.1]. 
\end{proof}
\begin{proposition} \label{specht fil dipdt}
Let $\mathrm{char}~K \neq 2,3.$ Then the permutation module $M^{\lambda,\mu}$ of $K\Sg_{a,b}$ has a dual Specht filtration.
\end{proposition}
\begin{proof}
If $\mathrm{char}~K \neq 2,3$, then by Lemma \ref{sym grp has dual Specht filtration}, we have $M^{\lambda} \in \mathcal{F}_{K\Sg_a}(S)$ and $M^{\mu} \in \mathcal{F}_{K\Sg_b}(S)$. From Lemma \ref{tensor product of two filtration for walled Brauer}, it follows that $M^{\lambda,\mu}=M^{\lambda}\boxtimes M^{\mu}\in \mathcal{F}_{K\Sg_a}(S) \times \mathcal{F}_{K\Sg_b}(S)$. This implies that the subquotient of $M^{\lambda,\mu}$ are of the form $S_{\lambda_i} \boxtimes S_{\mu_j}$ for some $i$ and $j$. By Lemma \ref{dual Specht module of direct product is tensor product of dual Specht module for direct product}, $S_{\lambda_i} \boxtimes S_{\mu_j}$ is isomorphic to the dual Specht module $S_{\lambda_i,\mu_j}$, for some $i,j$. Hence, $M^{\lambda,\mu}$ has a dual Specht filtration.
\end{proof}
The category of $K\Sg_{a,b}$-modules $X=X_1 \boxtimes X_2$ admitting a dual Specht filtration, $X_1 \in \mathcal{F}_{K\Sg_a}(S) $ and $X_2 \in \mathcal{F}_{KS_b}(S)$, is denoted by $\mathcal{F}_{K\Sg_{a,b}}(S)$.
Let $A$ is $K\Sg_{a,b}$-module. We say $A$ is \textit{relative projective} in $\mathcal{F}_{\Sg_{a,b}}(S)$ if $\text{Ext}_{K\mathfrak{S}_{a,b}}^{1} (A,X)=0$ for all $X \in \mathcal{F}_{K\mathfrak{S}_{a,b}}(S)$. 
\begin{proposition}\label{reproj}
If $\mathrm{char} ~K\neq 2,3$, and let $(\lambda, \mu) \in \Lambda_{a,b}$, then the permutation module $M^{\lambda,\mu}$ of $K\mathfrak{S}_{a,b}$ is relative projective in $\mathcal{F}_{K\mathfrak{S}_{a,b}}(S)$.
\end{proposition}
\begin{proof}
It suffices to check that $\text{Ext}_{K\mathfrak{S}_{a,b}}^{1} 
(M^{\lambda,\mu},X)=0$ for all $X=( X_1 \boxtimes X_2) \in \mathcal{F}_{K\mathfrak{S}_{a,b}}(S)$. By [\cite{CE_B}, Chapter XI, Theorem 3.1] we have 
\begin{align*}
\Ext_{K\mathfrak{S}_{a,b}}^{1} ( M^{\lambda,\mu}, X) &\cong \big(\Ext_{K\mathfrak{S}_{a}}^{1} (M^{\lambda}, X_1)\otimes_K \Hom_{K\mathfrak{S}_{b}} (M^{\mu},X_2) \big) \oplus \\
&\big( \Hom_{K\mathfrak{S}_{a}} (M^{\lambda},X_1)\otimes_K \Ext_{K\mathfrak{S}_{b}}^{1} (M^{\mu}, X_2)\big).
\end{align*}
Since the permutation module of the symmetric group is relative projective by [\cite{HN}, Proposition 4.1.1], we have $\Ext_{K\mathfrak{S}_{a}}^{1}(M^{\lambda},X_1)=0$ and $\Ext_{K\mathfrak{S}_{b}}^{1}(M^{\mu},X_2)=0$. Thus, each summand vanishes. 
\end{proof}

\begin{corollary}\label{relative projective for KSg}
 If $\mathrm{char}~K \neq 2,3$, then the Young modules are relative projective in $\mathcal{F}_{K\mathfrak{S}_{a,b}}(S)$.
\end{corollary}
\begin{proof}
This follows from the fact that the direct summands of relative projective modules are also relative projective. 
\end{proof}

Let $A$ be a relative projective module in $\mathcal{F}_{K\mathfrak{S}_{a,b}}(S)$. We call $A$ a \textit{relative projective cover} of $A'\in \mathcal{F}_{K\mathfrak{S}_{a,b}}(S)$ if it is minimal with the property that there exists an epimorphism 
\begin{tikzcd}
A \ar[r,twoheadrightarrow,"f"] & A'
\end{tikzcd}
with $\ker f\in \mathcal{F}_{K\mathfrak{S}_{a,b}}(S)$.
\begin{proposition}\label{relative projective cover for KSg}
Let $\mathrm{char}~K \neq 2,3$. Then the Young modules $Y^{\lambda,\mu}$ of $K\Sg_{a,b}$ are the relative projective covers of the dual Specht modules in $\mathcal{F}_{K\mathfrak{S}_{a,b}}(S)$. 
\end{proposition}
\begin{proof}
Recall that $Y^{\lambda,\mu}$ is the relative projective in $\mathcal{F}_{K\mathfrak{S}_{a,b}}(S)$ by Corollary \ref{relative projective for KSg}. We need to show that there exists a surjective homomorphism 
\begin{tikzcd}
Y^{\lambda,\mu} \ar[r,twoheadrightarrow,"f"] & S_{\lambda,\mu}
\end{tikzcd}
with $\ker f \in \mathcal{F}_{K\mathfrak{S}_{a,b}}(S)$.  Since $Y^{\lambda}$ and $Y^{\mu}$ are relative projective covers of $S_{\lambda}$ and $S_{\mu}$ respectively, we have surjective homomorphisms
\begin{tikzcd}
Y^{\lambda} \ar[r,twoheadrightarrow,"f_1"] & S_{\lambda},
\end{tikzcd}
and \begin{tikzcd}
Y^{\mu} \ar[r,twoheadrightarrow,"f_{2}"] & S_{\mu}
\end{tikzcd}
with $\ker f_1 \in \mathcal{F}_{\Sg_{a}}(S)$, and $\ker f_2 \in \mathcal{F}_{\Sg_{b}(S)}$, respectively. Thus, identifying $f$ with $f_1 \boxtimes f_2$,  which makes $f$ surjective. Furthermore, $\ker f= \ker f_1 \boxtimes \ker f_2$, which lies in $\mathcal{F}_{\Sg_{a}}(S) \times \mathcal{F}_{\Sg_{b}}(S)= \mathcal{F}_{K\mathfrak{S}_{a,b}}(S)$ by Lemma \ref{tensor product of two filtration for walled Brauer}. The minimality condition follows from the relative projectivity of $Y^{\lambda,\mu}$, and the Young modules being indecomposable.
\end{proof}

\section{Permutation module of $\B_{r,t}(\delta)$} \label{Permutation module for Walled Brauer algebra}
This section aims to construct the permutation and Young modules of $\B_{r,t}(\delta)$ using induction functors. The primary goal of this section is to establish both the existence and uniqueness of the Young modules, and it concludes with a characterization of these modules.

\subsection{Induction and Restriction functor}\label{Functors} 
In this section, we use an exact split pair of functors to compare the cohomology between two categories of modules. This approach was developed by Diracca and Koenig in \cite{DK}, and it also appears in the context of cleft extensions by Beligiannis in \cite{B2000}. Later on it is discussed in the case of cellularly stratified algebras in \cite{HHKP}. As shown in Section \ref{Walled Brauer algebras}, the walled Brauer algebras are cellularly stratified, so we can use this technique to define the induction and restriction functors. For $0 \leq l \leq s$, define 
\begin{align*}
\ind_l : \text{ \textbf{mod}-} K\mathfrak{S}_{r-l,t-l}  &\longrightarrow  \text{\textbf{mod}-} B \\
M &\longmapsto  M \otimes_{\substack{K\mathfrak{S}_{r-l,t-l}}} e_l(B/J_{l+1}),
\end{align*}
where $\text{\textbf{mod}-}K\mathfrak{S}_{r-l,t-l}$ and $ \text{\textbf{mod}-}B$ respresent the categories of right $K\mathfrak{S}_{r-l,t-l}$-modules and right $B$-modules, respectively. These functors allow us to compare the $\text{Ext}$ [\cite{HHKP}, Corollary 7.4] and the $\Hom$ [\cite{HHKP}, Proposition 4.3] of these modules. As the input algebra $K\mathfrak{S}_{r-l,t-l}$ is isomorphic to a subalgebra of $e_lBe_l$, this motivates us to define the exact localization and the right exact globalization functors. These are similar to the functors introduced by Green in \cite{Gr}. For $0 \leq l \leq s$, define
\begin{align*}
 \Ind_l : \text{ \textbf{mod}-}K\mathfrak{S}_{r-l,t-l}  &\longrightarrow  \text{\textbf{mod}-}B    &\Res_l : \text{\textbf{mod}-} B &\longrightarrow \text{ \textbf{mod}-}K\mathfrak{S}_{r-l,t-l}  \\
 M &\longmapsto  M \otimes_{K\mathfrak{S}_{r-l,t-l}} e_lB   &N &\longmapsto  N\otimes_{B} Be_l.
\end{align*}

We now define a poset $\Lambda$ as follows: $$\Lambda:= \{ (l,(\lambda,\mu)): 0 \leq l \leq s, \text{ and } (\lambda, \mu) \in \Lambda_{r-l,t-l} \}.$$ For any $(\lambda, \mu), (\lambda', \mu') \in \Lambda_{r-l,t-l}$, we can define a partial order $`` \leq "$ on $\Lambda$ as $$(l,(\lambda,\mu)) \leq (m,(\lambda',\mu')) \text{ if } m < l \text{ or } m=l, \text{ and } (\lambda, \mu) \unrhd (\lambda',\mu').$$

When $\delta =0$, then $\B_{1,1}(\delta)$ is not cellularly stratified. In this article, we consider all $B$-modules except for those when $B$ is not $\B_{1,1}(\delta)$ with $\delta=0$. 
\begin{proposition} [\cite{HHKP}, Proposition 4.2]\label{cell modules for B}
 Let $(l,(\lambda,\mu)) \in \Lambda$. The cell modules of $B$ are of the form $\ind_{l} S_{\lambda,\mu}$, where $S_{\lambda,\mu}$ are the dual Specht modules of $K\mathfrak{S}_{r-l,t-l}$. In particular, all simple $B$-modules are indexed by $(l,(\lambda,\mu))$, where $\lambda$ and $\mu$ are the $p$-regular partitions of $r-l$ and $t-l$,  respectively.
\end{proposition}
Let $(l,(\lambda,\mu)) \in \Lambda$ and $M^{\lambda,\mu}$ be the permutation module of $K\mathfrak{S}_{r-l,t-l}$. The \textit{permutation module} $M(l,(\lambda,\mu))$ of $B$ is defined as $M(l,(\lambda,\mu)):=\Ind_{l} M^{\lambda,\mu}$.

\subsection{Young modules}In this subsection, we define the Young module of the walled Brauer algebras, which is analogous to the Young module of the Brauer algebras as defined in \cite{HP}. Before describing the Young modules, we first present a theorem that establishes their existence and uniqueness within the context of $B$.

\begin{theorem}\label{exi young}
For $(l,(\lambda,\mu)) \in \Lambda$, there exists a unique direct summand of $M(l,(\lambda,\mu))$ whose quotient isomorphic to $\ind_l Y^{\lambda,\mu}$.
\end{theorem}
\begin{proof}
Using Theorem \ref{de per dirpdt}, we decompose the permutation module $M^{\lambda,\mu}$ of $K\mathfrak{S}_{r-l,t-l}$ into indecomposable Young modules. Applying $\Ind_l$ to this decomposition, we get $$\Ind_l M^{\lambda,\mu} = \bigoplus_{(l,(\lambda',\mu')) \in \Lambda} (\Ind_l Y^{\lambda',\mu'})^{\oplus a_{\lambda',\mu'}},$$ where $a_{\lambda',\mu'}$ is the multiplicity of $Y^{\lambda',\mu'}$ in $M^{\lambda,\mu}$. For $(l,(\lambda,\mu)) \in \Lambda$, we can further decompose $\Ind_l Y^{\lambda,\mu} = \bigoplus\limits_i (\Ind_l Y^{\lambda,\mu})\epsilon_i$, where $\epsilon_i$ are the primitive idempotents in $\End_{B}(\Ind_l Y^{\lambda,\mu})$ such that $1_{\End_{B}(\Ind_l Y^{\lambda,\mu})} = \sum\limits_i \epsilon_i$. The next step is to classify these indecomposable summands. To classify the indecomposable summands of $\Ind_l M^{\lambda,\mu}$, we need to answer the following questions:

\begin{itemize}
\item Investigate the direct summands of $\Ind_l Y^{\lambda,\mu}$, denoted by $(\Ind_lY^{\lambda,\mu}) \epsilon_i$, such that the quotient of this summand is isomorphic to $\ind_l Y^{\lambda,\mu}$.
\item Determine if the above direct summand is unique. 
\item Show that there are no other direct summands of $\Ind_l Y^{\lambda',\mu'}$ with a quotient that is isomorphic to $\ind_l Y^{\lambda,\mu}$, for $(\lambda',\mu')\neq(\lambda,\mu)$.
\end{itemize} 
The proof of the first two parts uses the same techniques as in Theorem 1 of \cite{In}; here we only provide a proof for the last part.\\
\textbf{Claim:} there does not exist any surjective map from $\Ind_l Y^{\lambda',\mu'}$ to $\ind_l Y^{\lambda,\mu}$ for $(\lambda',\mu')\neq (\lambda,\mu)$.

Assume there exists a direct summand of $\Ind_l Y^{\lambda',\mu'}$ such that we have a surjective map $\Phi : \Ind_l Y^{\lambda',\mu'} \longrightarrow \ind_l Y^{\lambda,\mu}$ for $(\lambda', \mu' ) \vartriangleright (\lambda,\mu) $. Since  
\begin{equation} \label{Hom between two Ind and ind}
\Hom_B ( \Ind_l Y^{\lambda',\mu'},\ind_l Y^{ \lambda,\mu} ) \cong \Hom_{K\mathfrak{S}_{r-l,t-l}} ( Y^{\lambda',\mu'}, Y^{\lambda,\mu}),
\end{equation}
 it follows that the induced map $\phi: Y^{\lambda',\mu'} \longrightarrow Y^{\lambda,\mu}$ is surjective with $(\lambda', \mu' ) \vartriangleright (\lambda,\mu) $. Consequently, we can extend $\phi$ linearly to $M^{\lambda', \mu'}$ in such a way that it takes all the indecomposable summands of $M^{\lambda', \mu'}$ to zero, except for $Y^{\lambda,\mu}$.

Let $k_{\lambda,\mu}$ be the alternating column sum of the $(\lambda,\mu)$-tableau $\mathfrak{t}_{\lambda,\mu}$, and let $e_{\lambda,\mu}$ be the polytabloid defined as in Subsection \ref{prop dipdt}. If $(\lambda', \mu') \triangleright (\lambda,\mu)$, then $(M^{\lambda'} \boxtimes M^{\mu'}) k_{\mathfrak{t}_{\lambda,\mu}} = (M^{\lambda'} k_{t_{\lambda}}) \boxtimes ( M^{\mu'}k_{t_{\mu}}) = 0$ by [\cite{JaB}, Lemma 4.6]. Then applying $\phi$, we have $\phi(M^{\lambda',\mu'} k_{\mathfrak{t}_{\lambda,\mu}})=0$. On the other hand, $ \phi(M^{\lambda',\mu'} k_{\mathfrak{t}_{\lambda,\mu}}) = \phi (M^{\lambda',\mu'}) k_{\mathfrak{t}_{\lambda,\mu}} = Y^{\lambda,\mu}k_{\mathfrak{t}_{\lambda,\mu}} $. But $Y^{\lambda,\mu}k_{\mathfrak{t}_{\lambda,\mu}}$ contains the generator $e_{\lambda,\mu}$ of $S^{\lambda,\mu}$. In particular, $Y^{\lambda,\mu}k_{\mathfrak{t}_{\lambda,\mu}}= (Y^{\lambda} k_{\mathfrak{t}_{\lambda}}) \boxtimes (Y^{\mu} k_{\mathfrak{t}_{\mu}}) \neq 0$ by [\cite{JaB}, Corollary 4.7], which is a contradiction. Thus, completing the proof. 
\end{proof}

\begin{definition}
Let us define the \textit{Young module} of $B$ to be the unique direct summand of $\Ind_l Y^{\lambda,\mu}$ with quotient isomorphic to $\ind_l Y^{\lambda,\mu}$. We denote this module by $Y(l,(\lambda,\mu))$.
\end{definition}

\begin{proposition} \label{Charaterize the Young modules for B}
The Young modules satisfy the following properties:
\begin{itemize}
\item[1.] The Young module $Y(l,(\lambda,\mu))$ appears exactly once as a direct summand in the decomposition of the permutation module $M(l, (\lambda,\mu))$ of $B$.
\item[2.] If $(l,(\lambda,\mu)) , (m,(\lambda',\mu')) \in \Lambda$ with $l>m$, then $Y(m,(\lambda',\mu'))$ does not appear as a summand of $M(l,(\lambda,\mu))$. 
\item[3.] If $(l,(\lambda,\mu)), (l,(\lambda',\mu')) \in \Lambda$, then  $Y(l,(\lambda',\mu'))$ occurs in the decomposition of $M(l,(\lambda,\mu))$ if and only if  $Y^{\lambda',\mu'}$ appears as a direct summand of $M^{\lambda',\mu'}$ which is possible only if $(\lambda', \mu' ) \unrhd (\lambda,\mu) $. Moreover, the multiplicity of $Y(l,(\lambda',\mu'))$ in $M(l,(\lambda,\mu))$ is equal to the multiplicity of $Y^{\lambda',\mu'}$ in $M^{\lambda,\mu}$.
\item[4.] $Y(l,(\lambda,\mu))\cong Y(m,(\lambda',\mu'))$ if and only if $(l,(\lambda, \mu))=(m,(\lambda',\mu'))$.
\end{itemize}
\end{proposition}
\begin{proof}
\begin{itemize}
\item[1.] This is an immediate consequence of Theorem \ref{exi young}.
\item[2.] Let $Y(m,(\lambda',\mu'))$ appears as a summand of $M(l,(\lambda,\mu))$ with $l>m$. There exists a surjective $B$-module homomorphism $ M(l,(\lambda,\mu)) \longrightarrow Y(m,(\lambda',\mu'))$, which can be extended to $ M(l,(\lambda,\mu)) \twoheadrightarrow \ind_{m} Y^{\lambda',\mu'} $. Therefore, we have
\begin{align*}
\Hom_{B} ( \Ind_{l} &M^{\lambda,\mu},  \ind_m  Y^{\lambda', \mu'}) \cong \Hom_{K\mathfrak{S}_{r-l,t-l}} (  M^{\lambda,\mu}, \Res_{l}\ind_m  Y^{\lambda', \mu'})\\
&\cong \Hom_{K\mathfrak{S}_{r-l,t-l}} (  M^{\lambda,\mu},  Y^{\lambda', \mu'}\otimes_{K\Sg_{r-m,t-m}}e_m (B/J_{m+1}) e_l) .
\end{align*}
If $l>m$, then $e_m (B/J_{m+1}) e_l  =0$. Hence, such a map does not exist.  

\item[3.] 
Let $Y(l,(\lambda',\mu'))$ occurs as a direct summand of $M(l,(\lambda,\mu))$. Then by an argument analogous to that in (\ref{Hom between two Ind and ind}), it follows that $Y^{\lambda',\mu'}$ is a direct summand of $M^{\lambda,\mu}$. By Theorem \ref{de per dirpdt}(1), we deduce that $(\lambda',\mu') \unrhd (\lambda,\mu)$. Conversely, if $Y^{\lambda',\mu'}$ is a direct summand of $M^{\lambda,\mu}$, then again by a similar argument to (\ref{Hom between two Ind and ind}), $Y(l,(\lambda',\mu'))$ is a direct summand of $M(l,(\lambda,\mu))$.

\item[4.]
Assume that $Y(l,(\lambda,\mu)) \cong Y(m, (\lambda',\mu'))$. Let $Y(l,(\lambda,\mu))$ is a direct summand of $M(l,(\lambda,\mu))$. Since $Y(l,(\lambda,\mu)) \cong Y(m, (\lambda',\mu'))$, it follows that $Y(m,(\lambda',\mu'))$ is also a direct summand of $M(l,(\lambda,\mu))$. By Proposition \ref{Charaterize the Young modules for B} (2), thich implies $m\geq l$. Similarly by Proposition \ref{Charaterize the Young modules for B} (3), if $l=m$, then $Y^{\lambda',\mu'}$ is a direct summand of $M^{\lambda,\mu}$ which implies $(\lambda',\mu') \unrhd (\lambda,\mu)$. Therefore, $(l,(\lambda,\mu)) \geq (m,(\lambda', \mu'))$. Reversing the roles,  we obtain $(l,(\lambda,\mu)) \leq(m,(\lambda',\mu'))$, which leads to the conclusion that $(l,(\lambda,\mu)) = (m,(\lambda',\mu'))$. 
\end{itemize}
\end{proof}

\section{Main Theorem} \label{imp sec}
In this section, we analyze the decomposition of the permutation modules of $B$ into indecomposable direct summands. We use a method similar to that used in \cite{In} and \cite{HHKP}. 

Let $\mathcal{F}_B(\Theta)$ \label{cell filtration defn} denote the category of $B$-modules that admit a cell filtration. By a \textit{cell filtration} we mean that if $M \in \mathcal{F}_B(\Theta)$, then it has a chain of submodules $0 \subseteq M_n \subseteq M_1 \subseteq \cdots \subseteq M_0=M$ such that each of the subquotient is isomorphic to a cell module. Note that $\ind_l S_{\lambda,\mu}$ are the cell modules of $B$ for $(l,(\lambda,\mu))\in \Lambda$, and we use $\Theta$ to denote the set of all cell modules of $B$.  

Let us recall the definition of a \textit{standard system} from [\cite{HHKP}, Definition 10.1]. 
\begin{definition} \label{definition of standard system}
Let $A$ be any algebra, and suppose we are given a finite set of non-isomorphic $A$-modules $\Theta(j)$, indexed by $j \in I$, where $I$ is equipped with a partial order $\leq$. Then the modules $\Theta(j)$ are said to form a standard system if the following
three conditions hold: 
\begin{itemize}
\item[(i)] For all $j \in I$, $\End_A (\Theta(j))$ is a division ring.
\item[(ii)] For all $m, n \in I$, if $\Hom_A(\Theta(m), \Theta(n)) \neq 0$, then $m \geq n$.
\item[(iii)] For all $m,n\in I$, if $\Ext_A^1(\Theta(m), \Theta(n)\neq 0$, then $m >n$.
\end{itemize}
\end{definition}

\begin{theorem}\label{stand sys}
If $\mathrm{char}~K \neq 2,3$, then the cell modules of $B$ form a standard system. 
\end{theorem}
\begin{proof}
We begin by proving the second condition for standard system, as defined in Definition \ref{definition of standard system}. By applying [\cite{Xi}, Lemma 3.2], we have
\begin{align}\label{Hom between two Specht modules}
\Hom_{K\mathfrak{S}_{r-l,t-l}} (S^{\lambda',\mu'}, S^{\lambda, \mu})
 \cong \Hom_{K\mathfrak{S}_{r-l}} (S^{\lambda'}, S^{\lambda}) \boxtimes_K \Hom_{K\mathfrak{S}_{t-l}} (S^{\mu'}, S^{\mu}). 
\end{align}
This implies $\Hom_{K\Sg_{r-l,t-l}}( S^{\lambda',\mu'}, S^{\lambda, \mu}) $ is zero if either tensor is zero. When $\mathrm{char}~ K \neq 2$, by [\cite{JaB}, Corollary 13.17] we know that whenever $\lambda' \ntrianglerighteq \lambda$ or $\mu' \ntrianglerighteq \mu$, then (\ref{Hom between two Specht modules}) becomes zero. 

Moving on to the third condition we use [\cite{CE_B}, Chapter XI, Theorem 3.1] to obtain the following:
\begin{align*}
\Ext_{K\mathfrak{S}_{r-l,t-l}}^{1} ( S^{\lambda',\mu'}, S^{\lambda, \mu}) &\cong \big(\Ext_{K\mathfrak{S}_{r-l}}^{1} (S^{\lambda'}, S^{\lambda})\otimes_K \Hom_{K\mathfrak{S}_{t-l}} (S^{\mu'}, S^{\mu}) \big) \\
&\quad\oplus \big( \Hom_{K\mathfrak{S}_{r-l}} (S^{\lambda'},S^{\lambda})\otimes_K \Ext_{K\mathfrak{S}_{t-l}}^{1} (S^{\mu'}, S^{\mu})\big).
\end{align*}
For $\text{Ext}_{K\mathfrak{S}_{r-l,t-l}}^{1} ( S^{\lambda',\mu'}, S^{\lambda, \mu})= 0$, both direct summands must be zero. If the first direct summand is zero. Then, atleast one of the tensors must be zero. If $\mathrm{char}~K \neq 2,3$, whenever $\mu' \ntrianglerighteq \mu$ then $\Hom_{K\mathfrak{S}_{t-l}} (S^{\mu'}, S^{\mu})=0$ by [\cite{JaB}, Corollary 13.17], or whenever $\lambda' \ntrianglerighteq \lambda $ we have $\Ext_{K\mathfrak{S}_{r-l}}^{1} (S^{\lambda'}, S^{\lambda})=0$ by [\cite{HN}, Proposition 4.2.1]. Similar conditions apply for the second direct summand to be zero. Hence, we conclude that $\Ext_{K\mathfrak{S}_{r-l,t-l}}^{1} ( S^{\lambda',\mu'}, S^{\lambda, \mu}) =0$, whenever $\lambda' \ntrianglerighteq \lambda$ or $\mu'\ntrianglerighteq \mu$.

The expression for $\End_{K\Sg_{r-l,t-l}} S^{\lambda,\mu}$ can be written as 
$$\End_{K\mathfrak{S}_{r-l,t-l}} S^{\lambda,\mu} \cong (\End_{K\mathfrak{S}_{r-l}} S^{\lambda})\boxtimes_K (\End_{K\mathfrak{S}_{t-l}} S^{\mu}).$$ By [\cite{JaB}, Corollary 13.17], each component is isomorphic to $K$, which implies that $\End_{K\mathfrak{S}_{r-l,t-l}} S^{\lambda,\mu}$ is a division ring.

Therefore, for each $l$, the set of Specht modules of $K\Sg_{r-l,t-l}$ form a standard system. Thus by [\cite{HHKP}, Section 10, Remark (4)], we have the dual Specht modules of $K\Sg_{r-l,t-l}$ form a standard system, for each $l$. In particular, the dual Specht modules are the cell modules of $K\Sg_{r-l,t-l}$. Hence, the cell modules of $B$ form a standard system by [\cite{HHKP}, Theorem 10.2].
\end{proof}

\subsection{Cell filtration} 
We will now prove that when $\mathrm{char}~K\neq 2,3$ the restriction of the cell module of $B$ to $K\Sg_{r-l,t-l}$ possesses a dual Specht filtration. 

\subsubsection{\textbf{Stabilizer of $v$ in $V_l$}} \label{stab pardiag}
For $0 \leq l \leq s$, let $v$ be an arbitrary element in $V_l$ with $l$ horizontal edges connecting the vertices from $\{ r-l+1, \cdots, r\}$ to $\{ r+1, \cdots, r+l\}$. Any $(\sigma,\tau) \in K\Sg_{l,l}$ acts from the right on the partial diagram $v$. Here $\sigma$ acts from the left side of the wall on the vertices $r-l+1, \cdots, r$ and $\tau$ acts on the vertices $r+1, \cdots, r+l$ from the right side of the wall, following the same concatenation rules as the action of $B$ on $v$. 

\begin{figure}[H]
\begin{center}
\begin{tikzpicture}[x=0.5cm,y=0.5cm]
\draw[-] (2,0).. controls(5.5,-1.25)..(9,0);
\draw[-] (3,0).. controls(5,-0.75)..(8,0);
\draw[-] (4,0).. controls(8,-1)..(12,0);
\draw[-] (5,0).. controls(8.5,-0.75)..(11,0);
\draw[-] (6,0).. controls(8.5,-0.25)..(10,0);
\draw[-] (7,0.5)--(7,-1.5);
\fill (1,0) circle[radius = 1pt];
\fill (13,0) circle[radius = 1pt];
\node at (-1,0) {$v=$};
\node at (1,0.5) {\tiny $1$};
\node at (2,0.5) {\tiny $2$};
\node at (3,0.5) {\tiny $3$};
\node at (4,0.5) {\tiny $4$};
\node at (5,0.5) {\tiny $5$};
\node at (6,0.5) {\tiny $6$};
\node at (8,0.5) {\tiny $7$};
\node at (9,0.5) {\tiny $8$};
\node at (10,0.5) {\tiny $9$};
\node at (11,0.5) {\tiny $10$};
\node at (12,0.5) {\tiny $11$};
\node at (13,0.5) {\tiny $12$};
\end{tikzpicture}
\begin{tikzpicture}[x=0.5cm,y=0.5cm]
\draw[-] (2,1.5)--(3,0);
\draw[-] (3,1.5)--(2,0);
\draw[-] (4,1.5)--(6,0);
\draw[-] (5,1.5)--(5,0);
\draw[-] (6,1.5)--(4,0);
\node at (0.5,1) {$=$};
\end{tikzpicture}
\caption{Identify a partial diagram $v$ as an element of $\Sg_l$.}\label{Identify a partial diagram $v_l$ as a element of S_l}
\end{center}
\end{figure}
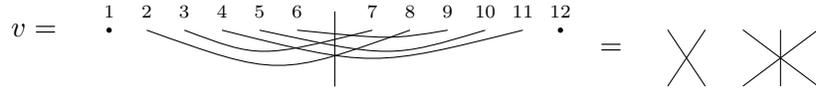

We want to determine the subgroup $ \Stab_{K\mathfrak{S}_{l,l}}(v)= \{ (\sigma,\tau) \in K\Sg_{l,l}: v.(\sigma,\tau)= v\}$. A partial diagram $v$ remains unchanged under two operation 
\begin{itemize}
\item[1.]Swapping one horizontal edge with another.
\item[2.] Exchanging the initial and terminal vertices of a horizontal edge.
\item[3.] Performing both operations simulteneously.
\end{itemize}
However, no element of $\Sg_{l,l}$ can be stabilize a partial diagram solely by interchanging the initial and terminal vertices of a horizontal edge. This is because such an operation requires the prmutation to belong to $\Sg_{2l}$ rather than $\Sg_{l,l}$. We now identify $v$ as a two-row diagram, which sends $v(i)$ to $j$ where $i$ and $j$ are the initial and terminal vertices of a horizontal edge, for $r-l+1 \leq i \leq r$, $r+1 \leq j\leq r+l$, for example see Figure \ref{Identify a partial diagram $v_l$ as a element of S_l}.

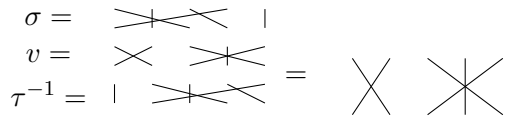
\begin{figure}[H]
\begin{center}
\begin{tikzpicture}[x=0.5cm,y=0.5cm]
\draw[-] (2,0)--(4,-0.5);
\draw[-] (3,0)--(3,-0.5);
\draw[-] (4,0)--(5,-0.5);
\draw[-] (5,0)--(2,-0.5);
\draw[-] (6,0)--(6,-0.5);
\draw[-] (2,-1)--(3,-1.5);
\draw[-] (3,-1)--(2,-1.5);
\draw[-] (4,-1)--(6,-1.5);
\draw[-] (5,-1)--(5,-1.5);
\draw[-] (6,-1)--(4,-1.5);
\draw[-] (2,-2)--(2,-2.5);
\draw[-] (3,-2)--(5,-2.5);
\draw[-] (4,-2)--(4,-2.5);
\draw[-] (5,-2)--(6,-2.5);
\draw[-] (6,-2)--(3,-2.5);
\node at (0.25,-1.25) {$v=$};
\node at (0.25,-0.25) {$\sigma=$};
\node at (0.25,-2.25) {$\tau^{-1}=$};
\end{tikzpicture}
\begin{tikzpicture}[x=0.5cm,y=0.5cm]
\draw[-] (2,1.5)--(3,0);
\draw[-] (3,1.5)--(2,0);
\draw[-] (4,1.5)--(6,0);
\draw[-] (5,1.5)--(5,0);
\draw[-] (6,1.5)--(4,0);
\node at (0.5,1) {$=$};
\end{tikzpicture}
\caption{An element of $\Stab_{K\Sg_{l,l}}(v).$}\label{An element of the stabilizer}
\end{center}
\end{figure}

Therefore, the action of $(\sigma,\tau)$ on $v$ (viewed as an element of $\mathfrak{S}_l$) is defined by $\sigma v\tau^{-1}$, which can be obtained by placing $\sigma$ at the top, followed by $v$, and $\tau^{-1}$ and concatenate them. Fixing $\sigma$ uniquely determines $\tau$. As any permutation in $\mathfrak{S}_l$ that swaps the horizontal edge in $v$ will fix $v$, therefore, $\Stab_{K\mathfrak{S}_{l,l}}(v)=K\Sg_{l}$. We denote the corresponding group algebra of this subgroup by $KH_l$. We illustrate this procedure with an example. Let $l=5$ and take $\sigma=(134) \in \Sg_5$. Then $\tau$ can be determined by $v^{-1}\sigma v$, so in this case $\tau =(254)$. So, $\sigma$ stabilizes $v$ by Figure \ref{An element of the stabilizer}.

\subsubsection{The structure of $J_me_l$}
For $0 \leq l< m \leq s$, let $e_l(J_m/J_{m+1})$ be a $K$-vector space consisting of walled Brauer diagrams with exactly $m$ horizontal edges that contain the partial diagram of $e_l$ at its top row. The left action of $\Sg_{r-l,t-l}$ on such diagram is as follows: $\Sg_{r-l}$ acts on the vertices $\{1,\cdots, r-l\}$ on the left side of the wall, and $\Sg_{t-l}$ acts on the vertices $\{r+l+1, \cdots, r+t\}$ on the right side of the wall. Therefore, $e_l(J_m/J_{m+1})$ is a left $\Sg_{r-l,t-l}$-module. 

\begin{lemma} \label{assum I}
As left $K\Sg_{r-l,t-l}$-modules, we have
$$ e_lJ_m \cong e_lJ_{m+1} \oplus  e_l(J_m/J_{m+1}) $$ for $0 \leq l< m \leq s$.
\end{lemma}
\begin{proof}
Consider the following filtration of the cell ideal $J_m$, for $0 \leq m \leq s$
\begin{equation} \label{exact sequence 1}
\begin{tikzcd}
0 \arrow[r] & J_{m+1} \arrow[r, "f"] & J_m \arrow[r, "g"] & J_m/ J_{m+1} \arrow[r] & 0,
\end{tikzcd}
\end{equation} 
where $f$ is an inclusion and $g$ is the quotient map. Hence, it is a short exact sequence of $B$-$B$ bimodules. Applying the exact functor $e_lB \otimes_{B} -$ on (\ref{exact sequence 1}), we obtain another short exact sequence 
 \begin{equation}\label{exact sequence 2}
 \begin{tikzcd}
0 \arrow[r] & e_lJ_{m+1} \arrow[r, "f' "] & e_lJ_m \arrow[r, "g'"] & e_l(J_m/ J_{m+1}) \arrow[r] & 0
\end{tikzcd}
\end{equation}
of $K\mathfrak{S}_{r-l,t-l}$-$B$ bimodules. As left $K\mathfrak{S}_{r-l,t-l}$ modules, we define an inclusion map $$g'': e_l(J_m/J_{m+1}) \xhookrightarrow{}  e_lJ_m.$$ Therefore, $g'\circ g''$ is an identity map on $e_l(J_m/J_{m+1})$, the sequence (\ref{exact sequence 2}) splits  as a sequence of $K\mathfrak{S}_{r-l,t-l}$-modules, and hence, the isomorphism is proved. 
\end{proof}
\textbf{The structure of $e_l (J_m/J_{m+1})e_m$}: Let $0 \leq l <m \leq s$. Consider $e_l(J_m/J_{m+1})e_m$ as a $K$-vector space spanned by all the walled Brauer diagrams that have exactly $m$ horizontal edges in each row with partial diagrams of $e_l$ and $e_m$ appearing in the top and bottom rows, respectively. The remaining $m-l$ horizontal edges in the top configuration can be arranged freely. The number of possible ways to arrange these $m-l$ horizontal edges is given by $|V_m^{m-l}|$, where $V_m^{m-l}$ is the subset of $V_m$ with the partial diagram of $e_l$ at the middle and $m-l$ horizontal edges arranged arbitrarily. The group algebra $K\Sg_{r-l,t-l}$ acts from the left on the vertices $\{1,\cdots, r-l\}$ and $\{r+l+1,\cdots, r+t\}$ in the top row, while $K\Sg_{r-m,t-m}$ acts from the right on the vertices $\{1,\cdots, r-m\}$ and $\{r+m+1,\cdots, r+t\}$ in the bottom row. Here,  $K\Sg_{r-l,t-l}$ permutes both the vertical edges as well the $m-l$ horizontal edges. Hence, the module $e_l(J_m/J_{m+1})e_m$ forms a free $K\Sg_{r-l,t-l}$-$K\mathfrak{S}_{r-m,t-m}$ bimodule with $K$-rank $|V_{m}^{m-l}|\cdot |K\Sg_{r-m,t-m}|$.

\begin{lemma} \label{assum II}
As left $K\mathfrak{S}_{r-l,t-l}$ modules, we have $$ e_l(J_m/J_{m+1}) \cong   e_l(B/J_{m+1})e_m \otimes_{e_mBe_m} e_m(B/J_{m+1})$$ for $0 \leq l < m\leq s$. 
\end{lemma}

\begin{proof}
The tensor product on the right is taken over $e_mBe_m$, which consists the walled Brauer diagram having at least $m$ horizontal edges. If $d \in e_m Be_m$ has more than $m$ horizontal edges, then its right action on $e_l(B/J_{m+1})e_m$ is zero because $e_l(B/J_{m+1})e_m$ consists only those diagrams with $m$ horizontal edges. An analogous statement holds for the left action of such a $d$ on $e_m(B/J_{m+1})$. Therefore, it suffices to consider only those diagrams in $e_mBe_m$ that have exactly $m$ horizontal edges. Thus, in this case $e_mBe_m$ is isomorphic to $K\Sg_{r-m,t-m}$. Hence, we may rewrite the tensor product as
\[
 e_l(B/J_{m+1})e_m \otimes_{e_mBe_m} e_m(B/J_{m+1}) \cong  e_l(B/J_{m+1})e_m \otimes_{K\Sg_{r-m,t-m}} e_m(B/J_{m+1}). 
\]
Moreover, as a $K\Sg_{r-l,t-l}$-$K\Sg_{r-m,t-m}$-bimodule, the module $e_l(B/J_{m+1})e_m$ is a free module with $K$-rank $|V_m^{m-l}|\cdot |K\Sg_{r-m,t-m}|$, that is $e_l(B/J_{m+1})e_m \cong K\Sg_{r-m,t-m}^{|V_m^{m-l}|}$. It follows that 
\begin{align}\label{dimesion are equal in assume II}
e_l(B/J_{m+1})e_m \otimes_{K\Sg_{r-m,t-m}} e_m(B/J_{m+1}) \cong \big( e_m(B/J_{m+1})\big)^{|V_m^{m-l}|}.
\end{align}
Define the $K$-linear map 
\begin{align}\label{surjective map in assume II}
e_l(B/J_{m+1})e_m \otimes_{K\Sg_{r-m,t-m}} e_m(B/J_{m+1}) &\longrightarrow e_l(J_m/J_{m+1})\nonumber \noindent\\
e_lde_m \otimes e_md' &\longmapsto e_lde_md',
\end{align}
where $d, d' \in B/J_{m+1}.$ Since $e_l(B/J_{m+1})$ consists of walled Brauer diagram with exactly $m$ horizontal edges containing the partial diagram of $e_l$ in the top, we have
\begin{align*}
\dim_K\big( e_l (J_m/J_{m+1})\big) & = |V_m^{m-l}|\cdot |V_m| \cdot |K \Sg_{r-m,t-m}|\\
&=\dim_K \big(e_m (B/J_{m+1})\big)^{|V_{m}^{m-l}|}.
\end{align*}
Since the map in (\ref{surjective map in assume II}) is surjective and the dimension of both the vector space are the same by (\ref{dimesion are equal in assume II}), it is an isomorphism of left $K\Sg_{r-l,t-l}$-modules.
\end{proof}

Let a group $G$ act on a set $\Sigma$, then the associated module $K[\Sigma]$ is called the $G$-permutation module. If $G$ acts transitively on $M$ and $H$ is the stabilizer of some element $m \in M$, then one can have $K[\Sigma]=K[H\backslash G]$, where $H\backslash G$ denotes the set of right cosets of $G$ modulo $H$.

Since $e_l (B/J_{m+1})e_m$ is a left $K\Sg_{r-l,t-l}$-module, then it follows that the left action of $K\Sg_{r-l,t-l}$ may4 fix the remaining $m-l$ horizontal edges. Therefore, we need to consider the stabilizer subgroup $KH_{m-l}$ of $K\Sg_{r-l,t-l}$. One can identify $K(\mathfrak{S}_{r-m,t-m} \times H_{m-l})$ as a subalgebra of $K\mathfrak{S}_{r-l,t-l}$, and is denoted by $H_{m,l}$. 

\begin{lemma} \label{L dec per} 
Let $0 \leq l <m \leq s$. As $K\mathfrak{S}_{r-l,t-l}$-$K\mathfrak{S}_{r-m,t-m}$-bimodules we have $$e_l(B/J_{m+1})e_m \cong K[\mathfrak{S}_{r-l,t-l}/H_{m-l}].$$
\end{lemma}

\begin{proof}
Consider the linear map 
\begin{equation*}
h: K[ \mathfrak{S}_{r-l,t-l}/H_{m-l}] \longrightarrow e_l(B/J_{m+1})e_m
\end{equation*}
given  by 
\begin{equation*}
(\sigma,\tau)H_{m-l} \longmapsto e_{m}^{\sigma,\tau},
\end{equation*}
where the walled Brauer diagram $e_{m}^{\sigma,\tau}$ is defined as follows: the top row of $e_{m}^{\sigma,\tau}$ contains the partial diagram of $e_l$ at its center and the bottom row contains the partial diagram of $e_m$ at its center. The remaining $m-l$ horizontal edges in the top row are connected from $\sigma(r-m+i)$ to $v^{-1}\tau v(r+m+1-i)$ for any $v\in V_m^{m-l}$, and $1 \leq i \leq m-l$. The remaining vertices in the top and bottom rows are connected with vertical edges by $\sigma(i)$ to $i$ for $1\leq i \leq r-m$, and $\tau(j)$ to $j$ for $r+m+1 \leq j \leq r+t$. This map is well-defined and surjective. Now, $\dim_{K}K[\Sg_{r-l,t-l}/H_{m-l}]=\frac{(r-l)!(t-l)!}{(m-l)!}$. On the other hand, we have  
\begin{align*}
\dim_K(e_m(B/J_{m+1})e_l)&=|V_m^{m-l}||\Sg_{r-m}||\Sg_{t-m}|\\
&= \binom{r-l}{m-l} \binom{t-l}{m-l} (m-l)! (r-m)! (t-m)!\\
&= \dim_K K[\Sg_{r-l,t-l}/H_{m-l}].
\end{align*}
Hence the map $h$ is an isomorphism of $K\Sg_{r-l,t-l}-K\Sg_{r-m,t-m}$-bimodules. 
\end{proof}

Let $\mathfrak{S}_{\lambda,\mu}$ denotes the Young subgroup of $\Sg_{r-l,t-l}$ for $(\lambda,\mu) \in \Lambda_{r-l,t-l}$. Since $H_{m,l}$ is a subalgebra of $K\Sg_{r-l,t-l}$, the set of double cosets $\{\mathfrak{S}_{\lambda,\mu} (\sigma,\tau) H_{m-l} : (\sigma,\tau) \in \Sg_{r-l,t-l}\}$. This forms an $\Sg_{r-m,t-m}$-set acting from the right, and this induces the $K\Sg_{r-m,t-m}$-permutation module, which is denoted by $K[\mathfrak{S}_{\lambda,\mu} \backslash \Sg_{r-l,t-l}/ H_{m-l} ]$. 
 The next lemma describes the $K\Sg_{r-m,t-m}$-orbits associated with a permutation module of $K\Sg_{r-m,t-m}$. Let $\mathfrak{T}$ be the collection of all compositions $(\gamma,\delta)$, such that $\gamma \vDash r-m$ and $\delta \vDash t-m$. Define $\pi:= \Sg_{\lambda,\mu}(\alpha,\beta)H_{m-l} \times \Sg_{r-m,t-m}$. 

\begin{lemma}\label{dec per 5}
For $m \geq l$, as $K\mathfrak{S}_{r-m,t-m}$-modules, we have $$K[\Sg_{\lambda,\mu}\backslash \mathfrak{S}_{r-l,t-l}/ H_{m-l} ]= \bigoplus_{\pi} M^{\phi(\pi)},$$ where $\phi$ is a function from $\mathfrak{S}_{r-m,t-m}$-orbits on $ \Sg_{\lambda,\mu}\backslash\mathfrak{S}_{r-l,t-l} \slash H_{m-l}$ to the set of compositions in $\mathfrak{T}$.
\end{lemma}
\begin{proof}
To show that the $\mathfrak{S}_{r-m,t-m}$-orbits on $  \Sg_{\lambda,\mu}\backslash\mathfrak{S}_{r-l,t-l} \slash H_{m-l}$ correspond to the permutation modules of $K\Sg_{r-m,t-m}$, it suffices to check that the stabilizer subgroups form a Young subgroup of $\Sg_{r-m,t-m}$.

Let $(\alpha,\beta) \in K\mathfrak{S}_{r-l,t-l}$, then we have
\begin{align*}
&\Stab_{K\mathfrak{S}_{r-m,t-m}}\big(  K\Sg_{\lambda,\mu}(\alpha,\beta)KH_{m-l}\big)\\
&= \{(\alpha',\beta') \in K\mathfrak{S}_{r-m,t-m}: \big( K\Sg_{\lambda,\mu} (\alpha,\beta)KH_{m-l}\big).(\alpha',\beta')=\big(K\Sg_{\lambda,\mu} (\alpha,\beta)KH_{m-l}\big)\}\\
&= \{(\alpha',\beta') \in K\mathfrak{S}_{r-m,t-m}: \big( K\Sg_{\lambda,\mu} (\alpha,\beta).(\alpha',\beta')KH_{m-l}\big)=\big(K\Sg_{\lambda,\mu} (\alpha,\beta)KH_{m-l}\big)\}\\
&= \Stab_{K\mathfrak{S}_{r-m,t-m}}\big(\Sg_{\lambda,\mu}(\alpha,\beta) \big)\\
&=K\mathfrak{S}_{r-m,t-m} \cap  \Stab_{K\mathfrak{S}_{r-l,t-l}}\big(K\Sg_{\lambda,\mu}(\alpha,\beta) \big )\\
&=K\mathfrak{S}_{r-m,t-m} \cap  (\alpha^{-1},\beta^{-1})
K\Sg_{\lambda,\mu}(\alpha,\beta)
\end{align*}
where the second equality holds because every element of $K\mathfrak{S}_{r-m,t-m}$ commute with $KH_{m-l}$. Define $\Sg_{\gamma,\delta}:=\mathfrak{S}_{r-m,t-m} \cap (\alpha^{-1},\beta^{-1})
\Sg_{\lambda,\mu}(\alpha,\beta)$, where $(\gamma,\delta) \in \mathfrak{T}$. Now, define a map $\phi: \Sg_{\lambda,\mu}\backslash\mathfrak{S}_{r-l,t-l} \slash H_{m-l} \times \Sg_{r-m,t-m} \longrightarrow \mathfrak{T}$ such that $\phi(\pi)= (\gamma,\delta)$. Hence, the permutation module corresponds to the double coset representative $(\alpha,\beta)$ is $M^{\gamma,\delta}$, which completes the proof.
\end{proof}

\subsubsection{Cell filtration of $M(l,(\lambda,\mu))$}
We have shown that if $\mathrm{char}~K \neq 2,3$, then $M^{\lambda,\mu}$ of $K\Sg_{a,b}$ has dual Specht filtration in Lemma \ref{specht fil dipdt}. This motivates us to study the cell filtration of the permutation module of $B$.
\begin{theorem}\label{cell fil}
If $\mathrm{char}~K \neq 2,3$, then the permutation module $M(l,(\lambda,\mu))$ of $B$ has a cell filtration.  
\end{theorem}
\begin{proof}
It is known that the algebra $B$ has a filtration by two-sided ideals $B = J_0 \supseteq J_{1} \supseteq \cdots \supseteq J_{s-1} \supseteq J_{s} \supseteq 0$. This filtration induces a short exact sequence of $B$-$B$ bimodules $0 \longrightarrow J_{m+1} \longrightarrow J_{m} \longrightarrow J_{m}/ J_{m+1} \longrightarrow 0$ for each $0\leq m \leq s$. Since the restriction functor $e_lB\otimes_B -$ is an exact functor, therefore, applying the functor to the above short exact sequence, we obtain a new short exact sequence $0 \longrightarrow e_lJ_{m+1} \longrightarrow e_lJ_{m} \longrightarrow e_l(J_{m}/ J_{m+1}) \longrightarrow 0$ of $e_l B e_l$-$B$ bimodules. By Lemma \ref{assum I}, the sequence is split exact as a left $K\mathfrak{S}_{r-l,t-l}$-module. Applying the left exact functor $ M^{\lambda,\mu} \otimes_{ K\mathfrak{S}_{r-l,t-l}} - $, we obtain the following exact sequence of right $B$-modules as follows:
\begin{align*}
0 \longrightarrow  M^{\lambda,\mu} \otimes_{ K\mathfrak{S}_{r-l,t-l}} e_l J_{m+1} \longrightarrow  &  M^{\lambda,\mu}\otimes_{ K\mathfrak{S}_{r-l,t-l}} e_lJ_{m}\\
&\longrightarrow  M^{\lambda,\mu} \otimes_{ K\mathfrak{S}_{r-l,t-l}}e_l(J_{m}/ J_{m+1}) \longrightarrow 0.
\end{align*}  
This leads to a filtration of $M(l,(\lambda,\mu)) = M^{\lambda,\mu}\otimes_{ K\mathfrak{S}_{r-l,t-l}} e_lB= \Ind_{l} M^{\lambda,\mu}$, given by
\begin{align*}
 M^{\lambda,\mu}\otimes_{ K\mathfrak{S}_{r-l,t-l}} e_lB \supseteq  M^{\lambda,\mu} &\otimes_{ K\mathfrak{S}_{r-l,t-l}} e_lJ_1 \supseteq \cdots \supseteq\\
&  M^{\lambda,\mu}\otimes_{ K\mathfrak{S}_{r-l,t-l}} e_lJ_{s-1}  \supseteq  M^{\lambda,\mu}\otimes_{ K\mathfrak{S}_{r-l,t-l}} e_lJ_s \supseteq 0
\end{align*}
 with the subquotient $ M^{\lambda,\mu} \otimes_{ K\mathfrak{S}_{r-l,t-l}}e_l(J_{m}/ J_{m+1}):= M^{m}(l,(\lambda,\mu))$ for $m=0,1,\cdots, s$. Now we observe that by Lemma \ref{assum II}, $ M^{m}(l,(\lambda,\mu)) =  M^{\lambda,\mu} \otimes_{ K\mathfrak{S}_{r-l,t-l}} \big( e_l(B/J_{m+1})e_m \otimes_{K\Sg_{r-m,t-m}} e_m(B/J_{m+1})\big)$. By using the Lemma~\ref{L dec per} on the term $e_l(B/J_{m+1})e_m$, and then applying $\ind_m$ we have $$\ind_{m} ( M^{\lambda,\mu} \otimes_{ K\mathfrak{S}_{r-l,t-l}}  K[\mathfrak{S}_{r-l,t-l} \slash H_{m-l}] ).$$ Using  the definition of $M^{\lambda,\mu}$ of $K\Sg_{r-l,t-l}$, we get 
\begin{align*}
& \ind_{m} ( (K[\Sg_{\lambda,\mu}\backslash \mathfrak{S}_{r-l,t-l}]) \otimes_{ K\mathfrak{S}_{r-l,t-l}} K[\mathfrak{S}_{r-l,t-l} \slash H_{m-l}]) \\
&= \ind_{m} \big(K[\Sg_{\lambda,\mu} \backslash \mathfrak{S}_{r-l,t-l}/H_{m-l} ] \big)\\
&= \ind_{m} \bigoplus_{\pi} M^{\phi(\pi)},
\end{align*} 
where the last equality follows from Lemma \ref{dec per 5}.
If $\mathrm{char}~K \neq 2,3$, then the Proposition \ref{specht fil dipdt} asserts that $M^{\phi(\pi)}$ has a dual Specht filtration. Since $\ind_m$ sends dual Specht modules to cell modules, we conclude that $M(l,(\lambda,\mu))$ has a cell filtration.
\end{proof}
\begin{corollary}
If $\mathrm{char}~K \neq 2,3$, then every direct summand of $M(l,(\lambda,\mu))$ also admits cell filtration. In particular, the Young module $Y(l,(\lambda,\mu))$ has a cell filtration.
\end{corollary}
\begin{proof}
By Theorem \ref{cell fil}, when $\mathrm{char}~K \neq 2,3 $, the module  $M(l,(\lambda,\mu))$ has a cell filtration and thus belongs to $\mathcal{F}_B(\Theta)$. We prove that the cell modules of $B$ form a standard system. Therefore, all direct summands (in particular Young module $Y(l,(\lambda, \mu))$) of $M(l,(\lambda,\mu))$ have cell filtration by [\cite{Ri}, Theorem 2].
\end{proof}

\subsubsection{Relative Projectivity}  
Relative projectivity is an important property of the permutation modules and the Young modules in the category of modules of the cellularly stratified algebras having with a cell filtration. Moreover, in the category of modules of quasi-hereditary algebras with a standard filtration, these modules are projective which will be useful in proving Theorem \ref{Dec Per}. Additionally, we show that the Young modules serve as a relative projective covers (as defined in \cite{HHKP}) of cell modules in $\mathcal{F}_B(\Theta)$.

\begin{lemma}\label{rela proj}
If $\mathrm{char}~K \neq 2,3$, then the restriction $\Res_l$ of a cell module of $B$ to $K\Sg_{r-l,t-l}$ is filtered by the dual Specht modules.
\end{lemma}

\begin{proof}
One can have $\Res_l \ind_n S_{\lambda, \mu}= S_{\lambda,\mu} \otimes_{K\Sg_{r-n,t-n}} e_n (B/J_{n+1})e_l$. If $n<l$, then $e_n (B/J_{n+1})e_l=0$. So, we assume $n \geq l$. For $n \geq l $ and $(\lambda,\mu) \in \Lambda_{r-n,t-n}$, we have 
\begin{align}
\Res_l \ind_n S_{\alpha,\beta} &\cong  S_{\lambda,\mu} \otimes_{K\Sg_{r-n,t-n}} e_n(B/J_{n+1})e_l\nonumber\\
&\cong S_{\lambda,\mu} \otimes_{K\Sg_{r-n,t-n}} K[H_{n-l}\backslash \Sg_{r-l,t-l}]~[\text{by Lemma }\ref{L dec per}]\nonumber\\
&\cong S_{\lambda,\mu} \otimes_{K\Sg_{r-n,t-n}} \Ind_{KH_{n-l}}^{K\Sg_{r-l,t-l}} 1_{H_{n-l}}\nonumber\\
&\cong S_{\lambda,\mu} \otimes_{K\Sg_{r-n,t-n}} \big( \Ind_{K\Sg_{n-l}}^{K\Sg_{r-l}} 1_{H_{n-l}} \boxtimes_K \Ind_{K\Sg_{1}}^{K\Sg_{t-l}} 1_{K\Sg_{1}}\big)\label{Ind commutes with tensors}\\
&\cong S_{\lambda,\mu} \otimes_{K\Sg_{r-n,t-n}} M^{(n-l,1^{r-n})}\boxtimes_K M^{(1^{t-l})}\nonumber\\
&\cong S_{\lambda,\mu} \otimes_{K\Sg_{r-n,t-n}} M^{(n-l,1^{r-n}),(1^{t-l})}\nonumber
\end{align} 
where (\ref{Ind commutes with tensors}) follows from [\cite{CR}, Lemma 10.17]. If $\mathrm{char}~K \neq 2,3$, the permutation module $M^{(n-l,1^{r-n}),(1^{t-l})}$ of $K\Sg_{r-n,t-n}$ has dual Specht filtrations by Proposition \ref{specht fil dipdt}. Hence, the lemma follows from the characteristic free version of Littlewood-Richardson rule . 
\end{proof}

\begin{theorem} \label{the permutation module is relative projective in B}
If $\mathrm{char}~K \neq 2,3$, then the permutation module is relative projective in $\mathcal{F}_B(\Theta)$. 
\end{theorem}

\begin{proof}
We need to show that $\Ext_{B} ( \Ind_{l} M^{\lambda,\mu}, C)=0 $ for all $C \in \mathcal{F}_B (\Theta)$.
    Consider the short exact sequence 
    \begin{equation} \label{short rela proj}
   \begin{tikzcd}
        0 \arrow[r] &C \arrow[r] & D \arrow[r, "\alpha"] &\Ind_l M^{\lambda,\mu} \arrow[r] &0
    \end{tikzcd}
\end{equation}   
   in $\Ext_{B}^{1} \big( M ((l,(\lambda,\mu)),C \big)$, where $C \in \mathcal{F}_B(\Theta)$. It is enough to show that (\ref{short rela proj}) splits. We first apply the exact functor $\Res_l$ in (\ref{short rela proj}),  followed by the left exact functor $\Hom_{K\Sg_{r-l,t-l}} (M^{\lambda,\mu}, -)$, then we get a long exact sequence 
\begin{equation}\label{long rela proj}
 \begin{tikzcd}
        0 \arrow[r]  &\Hom_{K\Sg_{r-l,t-l}} (M^{\lambda,\mu}, \Res_lC) \arrow[r]  &\Hom_{K\Sg_{r-l,t-l}} (M^{\lambda,\mu}, \Res_lD)\\
          \arrow[r]&\Hom_{K\Sg_{r-l,t-l}} (M^{\lambda,\mu}, \Res_l \Ind_{l}M^{\lambda,\mu}) \arrow[r] & 
        \Ext_{K\Sg_{r-l,t-l}}^{1} (M^{\lambda,\mu}, \Res_lC) \arrow[r] &\cdots.
    \end{tikzcd}
\end{equation}  
By Lemma \ref{rela proj}, since $C \in \mathcal{F}_B(\Theta)$, $\Res_l C$ has a dual Specht filtration, and by Corollary \ref{reproj} $M^{\lambda,\mu}$ is relative projective in $\mathcal{F}_{K\Sg_{r-l,t-l}} (S)$. It follows that $\Ext_{K\Sg_{r-l,t-l}}^{1}(M^{\lambda,\mu}, Ce_l)=0$. Hence, the long exact sequence (\ref{long rela proj}) terminates at $\Hom_{K\Sg_{r-l,t-l}} (M^{\lambda,\mu}, \Res_l\Ind_l M^{\lambda,\mu})$. By using the property that $\Res_l$ is right adjoint to $\Ind_l$, we obtain the following short exact sequence 
\begin{equation*}
  \begin{tikzcd}
          0 \arrow[r] & \Hom_{B} (\Ind_l M^{\lambda,\mu}, C) \arrow[r] & \Hom_{B} (\Ind_l M^{\lambda,\mu}, D)  \arrow[r, "\alpha'"] &\End_{B} (\Ind_l M^{\lambda,\mu}) \arrow[r] &0.
  \end{tikzcd}
\end{equation*}   
This implies $\alpha'$ is surjective. So, there exists $\beta: \Ind_l M^{\lambda,\mu} \longrightarrow D $ such that (\ref{short rela proj}) splits. 
\end{proof}

\begin{corollary}\label{proj cover}
If $\mathrm{char}~K \neq 2,3$, then all direct summands of the permutation module are relative projective in $\mathcal{F}_B(\Theta)$. Moreover, the Young module $Y(l,(\lambda,\mu))$ is the relative projective cover of the cell module $\ind_{l}S_{\lambda,\mu}$ in $\mathcal{F}_B(\Theta)$.
\end{corollary}
\begin{proof}
In the proof of Theorem \ref{the permutation module is relative projective in B}, replace $\Ind_l M^{\lambda,\mu}$ with the direct summand of $\Ind_l M^{\lambda,\mu}$ in (\ref{short rela proj}) we can then prove that each of the direct summand of $\Ind_l M^{\lambda,\mu}$ is relative projective in $\mathcal{F}_{B}(\Theta)$.

Next, we will prove that $Y(l,(\lambda,\mu))$ is relative projective cover of $\ind_{l}S_{\lambda,\mu}$. To show this we need to construct an epimorphism $\Psi: Y(l,(\lambda,\mu)) \longrightarrow \ind_{l}S_{\lambda,\mu}$ such that $\ker \Psi \in \mathcal{F}_{B}( \Theta)$ and $Y(l,(\lambda,\mu))$ is minimal with this property. 

Using the Proposition \ref{relative projective cover for KSg}, $S_{\lambda,\mu}$ is the top quotient of $Y^{\lambda,\mu}$ with the kernel lies in $\mathcal{F}_{K\Sg_{r-l,t-l}}(S)$. Since $\ind_l$ is an exact functor, the kernel of the epimorphism $\psi: \ind_l Y^{\lambda,\mu} \longrightarrow \ind_l S_{\lambda,\mu}$ belongs to $\mathcal{F}_B(\Theta)$. Note that the map $\psi$ is surjective by Proposition \ref{relative projective cover for KSg}, and the exactness of $\ind_l$  functor. Consider the following epimorphism $\phi \iota$ from Theorem \ref{exi young}
\begin{equation*}
\begin{tikzcd}
\phi  \iota: Y(l,(\lambda,\mu)) \arrow[r,"\iota"] &\Ind_l Y^{\lambda,\mu} \arrow[r, "\Phi"] &\ind_l Y^{\lambda,\mu}. 
\end{tikzcd}
\end{equation*}
where $\Phi$ and $\iota$ is surjective due to which $\Phi \iota$ is surjective.  The map $\Psi$ can be obtained by 
\begin{equation*}
\begin{tikzcd}
\Psi: Y(l,(\lambda,\mu)) \arrow[r,"\Phi\iota"] &\ind_l Y^{\lambda,\mu} \arrow[r, "\psi"] &\ind_l S_{\lambda,\mu}. 
\end{tikzcd}
\end{equation*} 
Since $\Phi\iota$ and $\psi$ both are surjective. Then so is $\Psi$. Consider the diagram 
\begin{center}
\begin{equation}\label{commutative diag for relative proj cover}
\begin{tikzcd}[contains/.style = {draw=none,"\in" description,sloped}]                                    
&0 &0\\
&\ker f \ar[r,dashrightarrow," \thicksim"]\ar[u,leftarrow] & \ker (\Phi\iota)\ar[u,leftarrow] \\
  0\ar[r] & \ker \Psi \ar[d,dashrightarrow,"f"]\ar[u,leftarrow] \ar[r] & Y(l,(\lambda,\mu)) \ar[d,"\Phi\iota"] \ar[r," \Psi"]\ar[u,leftarrow] & \ind_{l}S_{\lambda,\mu} \ar[d,equal] \ar[r] & 0 \\
  0\ar[r] & \ker \psi \ar[r]\ar[d] & \ind_l Y^{\lambda,\mu} \ar[r,"\psi"]\ar[d] & \ind_{l}S_{\lambda,\mu} \ar[r]        & 0 .\\
  & 0 &0 
\end{tikzcd}
\end{equation}
\end{center}
 Since the composition 
 \begin{equation*}
 \begin{tikzcd}
 \ker \Psi \arrow[r] &Y(l,(\lambda,\mu)) \arrow[r,"\Phi  \iota"] & \ind_l Y^{\lambda,\mu} \arrow[r, "\psi"] & \ind_{l}S_{\lambda,\mu}
 \end{tikzcd}
 \end{equation*}
is zero, and hence by the universal property of $\ker \psi$ there exists a unique morphism $f:\ker \Psi \longrightarrow \ker \psi$ such that (\ref{commutative diag for relative proj cover}) is commutative. Again, by the universal property of $\ker (\Phi \iota) $ and the snake lemma, we conclude that the map $\ker f \longrightarrow \ker (\Phi \iota)$ is an isomorphism. The snake lemma also implies that the map $f:\ker \Psi \longrightarrow \ker \psi$ is surjective. It follows that 
\begin{equation*}
\begin{tikzcd}
0 \arrow[r] & \ker (\Phi  \iota) \arrow[r] &\ker \Psi \arrow[r] &\ker \psi \arrow[r] &0 
\end{tikzcd}
\end{equation*}
is a short exact sequence with $\ker \psi \in \mathcal{F}_{B} (\Theta)$. We need to claim that $\ker \Phi  \iota \in \mathcal{F}_B(\Theta)$ which will then implies $\ker \Psi \in \mathcal{F}_B(\Theta)$ because $\mathcal{F}_B(\Theta)$ is extension-closed. Consider the commutative diagram

\begin{equation}\label{commutative diagram for extension closed portion}
\begin{tikzcd}[contains/.style = {draw=none,"\in" description,sloped}]                                    
  0\ar[r] & \ker (\Phi \iota) \ar[d,hookrightarrow,"\iota"] \ar[r] & Y(l,(\lambda,\mu)) \ar[d,hookrightarrow,"\iota"] \ar[r," \phi\iota"] & \ind_l Y^{\lambda,\mu} \ar[d,equal] \ar[r] & 0 \\
  0\ar[r] & \ker \Phi \ar[r] \ar[u, shift left," \pi"]& \Ind_l Y^{\lambda,\mu} \ar[r,"\phi"] \ar[u, shift left,"\pi"] & \ind_l Y^{\lambda,\mu} \ar[r]        & 0. \\
\end{tikzcd}
\end{equation}

Since $\iota \big( \ker (\Phi \iota)\big) \subset \ker \Phi$, this ensures that $\iota$ restricts a map $\ker(\Phi \iota) \hookrightarrow \ker \Phi$ in (\ref{commutative diagram for extension closed portion}). Also $Y(l,(\lambda,\mu))$ is a direct summand of $\Ind_l Y^{\lambda,\mu}$ by Theorem \ref{exi young}. So, we have a projection $\pi: \Ind_l Y^{\lambda,\mu} \longrightarrow Y(l,(\lambda,\mu))$ such that $\pi \big(\ker \Phi \big) \subset  \ker (\Phi \iota)$. Then $\pi$ restrict to a map $\ker \phi \longrightarrow \ker (\phi \iota)$. Thus, $\ker (\phi \iota)$ is a direct summand of $\ker \phi$. Since $\ker \Phi \in \mathcal{F}_B(\Theta) $, and $  \mathcal{F}_B(\Theta) $ is closed under the direct summands, it follows that $\ker (\Phi \iota) \in  \mathcal{F}_B(\Theta) $.

The minimality condition follows immediately because $Y(l,(\lambda,\mu))$ is relative projective in $\mathcal{F}_B(\Theta)$, and is indecomposable. This ends the proof. 
\end{proof}

\subsection{\textbf{Proof of the main Theorem }} The main theorem is proved using the standardization theorem, originally introduced by Dlab and Ringel in \cite{DR}. We will use the construction of the Schur algebra of a cellularly stratified algebra from \cite{HHKP} to achieve this.


\begin{proof}[Proof of Theorem~{\upshape\ref{Dec Per}}]
Since $B$ is a cellularly stratified algebra, then by [\cite{HHKP}, Theorem 13.1], there exists a quasi hereditary algebra $S(B):= \End_{B} Y$, where $Y= \oplus_{(l,(\lambda,\mu))\in \Lambda} M)L,(\lambda,\mu)),$ such that the following holds: the category $\mathcal{F}_{B}(\Theta)$ is equivalent to the category of $S(B)$-modules with a standard filtration, as exact categories. In \cite{DR}, Dlab and Ringel established an one-to-one correspondence between the modules in the standard system and the indecomposable relative projective modules in $\mathcal{F}_{B}(\Theta)$. By Theorem \ref{exi young}, for each $(l,(\lambda,\mu)) \in \Lambda$, there is a unique Young module, and a cell module indexed by $(l,(\lambda,\mu))$. By Proposition \ref{cell modules for B} and Theorem \ref{stand sys}, the cell modules are indexed by $\Lambda$ and form a standard system. From Corollary \ref{proj cover}, the Young module and all the direct summands are relative projective in $\mathcal{F}_B(\Theta)$. Using the one-to-one correspondence in \cite{DR}, we conclude that all the relative projective $B$-modules are indexed by $\Lambda$. The image of these relative projective $B$-modules under the equivalence corresponds to projective modules in the quasi-hereditary algebra $S(B)$ as stated in [\cite{DR}, Corollary 2]. Hence, the projective $S(B)$-modules are parametrized by elements of $\Lambda$. 

Since the permutation module $M(l,(\lambda,\mu))$ is relative projective in $\mathcal{F}_B(\Theta)$, its image  say $\mathcal{M}$ is a projective $S(B)$-module under the equivalence. Since $\mathcal{M}$ is projective, we can write $\mathcal M \cong \bigoplus\limits_{(n, ( \gamma,\nu))} \mathcal{Y}^{(n, ( \gamma,\nu))^{\oplus a_{(n, ( \gamma,\nu))}}}$ where $(n,(\gamma,\nu)) \in \Lambda$. So the image of $\mathcal{Y}^{(n, ( \gamma,\nu))}$ is relative projective in $\mathcal{F}_{B} (\Theta)$. By Corollary \ref{proj cover}, all the direct summands of $M(l,(\lambda,\mu))$ are relative projective in $\mathcal{F}_{B}(\Theta)$. Therefore, all $Y(m,(\lambda',\mu'))$ with $(m,(\lambda',\mu')) \in \Lambda$ appear in the decomposition of $M(l,(\lambda,\mu))$, but some of them may not occur, depending on the condition in Proposition \ref{Charaterize the Young modules for B} which is $(m,(\lambda',\mu')) \leq (l,(\lambda,\mu))$. 
\end{proof}

\begin{remark}
 The new Schur algebra of the walled Brauer algebra $\B_{r,t}(\delta)$ is defined as $$S(\mathcal{B}_{r,t}(\delta)):= \End_{\mathcal{B}_{r,t}(\delta)} \big (\bigoplus_{(l,(\lambda,\mu)) \in \Lambda }M(l,(\lambda,\mu))\big)$$ where $M(l,(\lambda,\mu))$ is the permutation module of $\B_{r,t}(\delta)$ analogous to those defined in \cite{HK2012}. Although the rational Schur algebra which is the Schur algebra of $\B_{r,t}(\delta)$ has been previously examined as a centralizer algebra acting on mixed tensor space in \cite{DD}, studying the new Schur algebra is still a compelling topic. In a forthcoming paper, we shall investigate the structure and properties of $S(\B_{r,t}(\delta))$.
\end{remark}

\textbf{Acknowledgements}

The authors express their gratitude to Prof. Steffen Koenig for his time, valuable discussions and proofreading of this article. The first author's research is supported by Indian Institute of Science Education and Research Thiruvananthapuram PhD fellowship. The second author's research was partially supported by IISER-Thiruvananthapuram, SERB-Power Grant SPG/2021/004200 and Prof. Steffen Koenig's research grant. The second author also would like to acknowledge the Alexander Von Humboldt Foundation for their support.


\begin{thebibliography}{0}
\bibitem{B2000} Beligiannis, A. Cleft extensions of Abelian categories and applications to ring theory, Comm. Algebra { \bf 28}, 4503–4546 (2000).
\bibitem{BCHL} Benkart, G., Chakrabarti, M., Halverson, T., Leduc, R., Lee, C., Stroomer, J. Tensor product representations of general linear groups and their connections with Brauer algebras. J. Algebra {\bf 166}(3), 529–567 (1994). 

\bibitem{BS} Brundan, J., Stroppel, C. Gradings on walled Brauer algebras and Khovanov’s
arc algebra. Adv. Math. {\bf 231}(2), 709–773 (2012). 
\bibitem{CE_B} Cartan, H., Eilenberg, S. Homological algebra (PMS-19), Volume 19. Prince-
ton University Press, Princeton (1956). 

\bibitem{CT} Chuang, J.,  Tan, K. M.,  Representations of wreath products of algebras, Math. Proc. Cambridge Philos. Soc. {\bf 135}(3), 395--411 (2003).
\bibitem{CVDM} Cox, A., De Visscher, M., Doty, S., Martin, P. On the blocks of the walled Brauer
algebra. J. Algebra 320(1), 169–212 (2008). 
\bibitem{CR} Curtis, C.W., Reiner, I. Methods of representation theory with applications to finite groups and orders. Vol. I. Pure and
Applied Mathematics, p. 819. John Wiley \& Sons, Inc., New York, (1981).
\bibitem{DD} Dipper, R., Doty, S. The rational Schur algebra. Represent. Theory {\bf 12}, 58–82
(2008). 
\bibitem{DK} Diracca, L., Koenig, S. Cohomological reduction by split pairs. J. Pure Appl.
Algebra {\bf 212}(3), 471–485 (2008).
\bibitem{DR} Dlab, V., Ringel, C.M. The module theoretical approach to quasi-hereditary algebras. In: Representations of Algebras and Related Topics (Kyoto, 1990). London
Math. Soc. Lecture Note Ser., vol. 168, pp. 200–224. Cambridge Univ. Press,
Cambridge, (1992).
\bibitem{Er} Erdmann, K. Young modules for symmetric groups. vol. 71, pp. 201–210 (2001).

\bibitem{Grb} Grabmeier, J. Unzerlegbare Moduln mit trivialer Youngquelle und Darstellungs-
theorie der Schuralgebra. Bayreuth. Math. Schr. (20), 9–152 (1985).
\bibitem{Gr} Green, J.A. Polynomial representations of $\textrm{GL}_n$ . Lecture Notes in Mathematics, vol. 830, p. 118. Springer-Verlag, Berlin–New York (1980).
\bibitem{Ru} Green, R. Some properties of specht modules for the wreath product of symmetric groups. PhD thesis, University of Kent, (May 2019). 
\bibitem{GM} Green, R.M., Martin, P.P. Constructing cell data for diagram algebras. J. Pure Appl. Algebra { \bf 209}(2), 551–569 (2007). 
\bibitem{HHKP} Hartmann, R., Henke, A., Koenig, S., Paget, R. Cohomological stratification of diagram algebras. Math. Ann. { \bf 347}(4), 765–804 (2010).
\bibitem{HP} Hartmann, R., Paget, R. Young modules and filtration multiplicities for
Brauer algebras. Math. Z. {\bf 254}(2), 333–357 (2006). 
\bibitem{H} Hemmer, D.J. Symmetric group modules with Specht and dual Specht fil-
trations. Comm. Algebra {\bf 35}(11), 3292–3306 (2007). 
\bibitem{HN} Hemmer, D.J., Nakano, D.K. Specht filtrations for Hecke algebras of type
A. J. London Math. Soc. (2) {\bf 69}(3), 623–638 (2004).
\bibitem{HK2012} Henke, A., Koenig, S. Schur algebras of Brauer algebras I. Math. Z. { \bf 272}(3-4),729–759 (2014). 
\bibitem{JaB} James, G.D. The representation theory of the symmetric groups. Lecture Notes in Mathematics, vol. 682, p. 156. Springer, Berlin, (1978).
\bibitem{Ja} James, G.D. Trivial source modules for symmetric groups. Arch. Math. (Basel)
{ \bf 41}(4), 294–300 (1983).
\bibitem{Kl} Klyachko, A.A. Direct summands of permutation modules. Selected translations. Selecta Math. Soviet. {\bf 3}(1), 45–55 (1983/84).
\bibitem{K} Koike, K. On the decomposition of tensor products of the representations of the classical groups: by means of the universal characters. Adv. Math. {\bf 74}(1), 57–86
(1989). 
\bibitem{KX99} Koenig, S., XI, C. Cellular algebras: Inflations and morita equivalences. J. London Math. Soc. { \bf 60}(3), 700–722 (1999). 
\bibitem{In} Paul, I. Permutation modules for cellularly stratified algebras. J. Pure Appl.
Algebra { \bf 224}(11), 106412–33 (2020). 
\bibitem{Ri} Ringel, C.M. The category of modules with good filtrations over a quasi-
hereditary algebra has almost split sequences. Math. Z. { \bf 208}(2), 209–223 (1991).

\bibitem{RS1} Rui, H., Su, Y. Affine walled Brauer algebras and super Schur-Weyl duality. Adv.
Math. { \bf 285}, 28–71 (2015).
\bibitem{T} Turaev, V.G. Operator invariants of tangles, and R-matrices. Izv. Akad.
Nauk SSSR Ser. Mat. { \bf 53}(5), 1073–1107, 1135; translation in
Math. USSR-Izv.35, no.2, 411–444 (1990). 
\bibitem{Xi} Xi, C. On the representation dimension of finite dimensional algebras. J. Algebra {\bf 226}(1), 332–346 (2000).

\end{thebibliography}

\end{document}